\documentclass[11pt,twoside]{article}
\usepackage{amssymb,amsfonts,amsmath,latexsym,graphicx}
\newtheorem{theorem}{Theorem}[section]
\newtheorem{lemma}[theorem]{Lemma}

\newtheorem{coro}[theorem]{Corollary}
\newtheorem{defin}[theorem]{Definition}

\renewcommand{\appendix}[2]{\section*{#1}\label{#2}\addtocontents{toc}{\vskip 0pt {\hspace*{-12pt}\bf #1\hfill\ \thepage}}}
\newcommand{\grad}{\nabla}
\newcommand{\xO}{\Omega}
\newcommand{\email}[1]{{\small E-mail: {\textsf {#1}}}}

\newcommand{\affil}[1]{{\small\sl #1}}

\newcommand\runninghead[1] {\pagestyle{myheadings}\markboth {{\footnotesize\it{\quad #1}\hfill}}{{\footnotesize\it{#1\hfill\quad}}}}\headsep=40pt

\newcommand{\be}{\begin{equation}}
\newcommand{\ee}{\end{equation}}

\newcommand{\RN}{\mathbb{R}^N}

\newcommand{\ana}{\nabla}
\newcommand{\CO}{C_{0}^{\infty} (\Omega)}

\newcommand{\CRN}{C_{0}^{\infty} (\RN)}

\newcommand{\bea}{\begin{eqnarray}}
\newcommand{\eea}{\end{eqnarray}}

\newcommand{\xD}{\Delta}

\newcommand{\xL}{\Lambda}

\newcommand{\xo}{\omega}
\newcommand{\XO}{\Omega}

\newcommand{\ia}{({\rm i})}
\newcommand{\ib}{({\rm ii})}

\renewcommand\runninghead[2]{\pagestyle{myheadings}\markboth
{{\footnotesize\it{\quad #1}\hfill}}
{{\footnotesize\it{#2\hfill\quad}}}}\headsep=40pt

\begin{document}
\runninghead{\hfill\today\hfill ---\hfill  A. Tertikas \& N. B.
Zographopoulos}{Best Hardy-Rellich inequalities}

\title{ Best Constants in the Hardy-Rellich Inequalities and Related Improvements }
\author{
{\bf A. Tertikas}\thanks{A.T. acknowledges partial support by the
RTN European network Fronts--Singularities,
HPRN-CT-2002-00274}\\
\affil{Department of Mathematics, University of Crete, 71409 Heraklion, Greece}\\
\affil{and Institute of Applied and Computational Mathematics,}\\
\affil{FORTH, 71110 Heraklion, Greece}\\
\email{tertikas@math.uoc.gr}\\
\\
{\bf N. B. Zographopoulos}\thanks{N.Z. acknowledges partial
financial support from the PYTHAGORAS Program No. 68/831 of the
Ministry of Education of the Hellenic Republic.
Part of this work was announced in \cite{tz03}.}\\
\affil{Department of Mathematics, University of Aegean,
83200 Karlovassi, Samos, Greece} \\
\email{zographopoulos@aegean.gr}\\
}
\date{}
\maketitle \thispagestyle{empty}
\begin{abstract}
We consider Hardy-Rellich inequalities and discuss their possible
improvement. The procedure is based on decomposition into
spherical harmonics, where in addition various new inequalities
are obtained (e.g. Rellich-Sobolev inequalities). We discuss also
the optimality of these inequalities in the sense that we
establish (in most cases) that the constants appearing there are
the best ones. Next, we investigate the polyharmonic operator
(Rellich and Higher Order Rellich inequalities); the difficulties
arising in this case come from the fact that (generally)
minimizing sequences are no longer expected to consist of radial
functions. Finally, the successively use of the Rellich
inequalities lead to various new Higher Order Rellich
inequalities.
\end{abstract}
\section{Introduction}
Hardy inequality states that for $N \geq 3,$ for all $ u \in
\CRN$\be \label{1.1}
 \int_{\RN} |\ana
u|^2 \, dx \geq \biggl (\frac{N-2}{2}\biggr )^2\, \int_{\RN}
\frac{u^2}{|x|^2} \, dx.
 \ee
 The constant $(\frac{N-2}{2})^2$ is
the best constant in inequality (\ref{1.1}). A similar inequality
with the same best constant
 holds if $\RN$   is replaced by $\XO$ and $\XO$
contains the origin.

When  $\XO$ is a bounded domain, a much stronger inequality was
discovered by Brezis and V\'{a}zquez \cite{bv}, that is for all $
u \in \CO,$ \be \label{1.2} \int_{\XO} |\nabla u|^2 dx \geq
\left(\frac{N-2}{2} \right)^2 \int_{\XO} \frac{u^2}{|x|^2} dx +
\,z_0^2  \, \left(\frac{\xo_N}{|\XO|} \right)^{\frac{2}{N}}
\int_{\XO} u^2 dx, \ee
 where $\xo_N$  and $|\XO|$
denote the volume of the unit ball and $\xO$ respectively, and
$z_0 = 2.4048 \ldots$ denotes the first zero of the Bessel
function $J_0(z)$. Inequality (\ref{1.2}) is optimal in case $\XO$
is a ball centered at zero. We set $D =  \sup_{x \in \XO} |x|$ and
define recursively
\begin{eqnarray}
&& X_1(t) = (1- \log t)^{-1},\quad t\in (0,1],\nonumber \\
&& X_k(t) = X_{1}(X_{k-1}(t)), ~~~~~~~k=2,3,\ldots\, , t\in (0,1].
\label{1.3}
\end{eqnarray}
In \cite{ft01} actually, the following improved Hardy inequality
was also established for $ u \in \CO$, \be \int_{\XO} |\nabla u|^2
dx \;\geq  \; \left(\frac{N-2}{2} \right)^2
 \int_{\XO} \frac{u^2}{|x|^2}  dx
 +  \frac{1}{4}  \sum_{i=1}^{\infty} \int_{\XO} \frac{u^2}{|x|^2}
 X_1^2 X_2^2
 \ldots X_i^2     dx,
\label{1.4} \ee where we use the notation $X_i$ for
$X_i(\frac{|x|}{D})$. We will make use of the same notation
throughout this work.   Here the constants that appear are best
constants. It is worth mentioning that for $N\geq 2m+2 \, > \, 2$
and $ u \in \CO$ inequality (\ref{1.4}) takes the equivalent form
\be \label{1.5} \int_{\XO} \frac{|\nabla u|^2}{|x|^{2m}} dx \;
\geq  \;\left(\frac{N-2m-2}{2} \right)^2 \int_{\XO}
\frac{u^2}{|x|^{2m+2}}  dx\; + \;  \frac{1}{4} \sum_{i=1}^{\infty}
\int_{\xO} \frac{u^2}{|x|^{2m+2}} X_1^2 X_2^2 \ldots X_i^2
 dx. \ee

Similarly to (\ref{1.1}), the classical Rellich inequality states
that for $N \geq 5,$ for all $ u \in \CRN$,
\be \label{rel1} \int_{\RN} (\xD u)^2 \, dx \geq \biggl
(\frac{N(N-4)}{4}\biggr )^2\, \int_{\RN} \frac{u^2}{|x|^4} \, dx.
\ee
Davies and Hinz \cite{dav98} obtained various Rellich inequalities
as well as higher order Rellich inequalities. Gazzola,  Grunau and
Mitidieri \cite{ggm03} on the other hand obtained improved Rellich
inequalities in the spirit of \cite{bv}. As an example we mention
the following inequality that holds true for $N \geq 5,$ and all $
u \in \CO$ \be \label{1.7} \int_{\XO} (\xD u)^2 dx \geq
\left(\frac{N(N-4)}{4} \right)^2 \int_{\XO} \frac{u^2}{|x|^4} dx +
\, \frac{N(N-4)}{2} \, \xL_2 \, \left(\frac{\xo_N}{|\XO|}
\right)^{\frac{2}{N}} \int_{\XO} \frac{u^2}{|x|^2} dx \, + \,
\xL_4^2 \, \left(\frac{\xo_N}{|\XO|} \right)^{\frac{4}{N}}
\int_{\XO} u^2 dx .\ee  Constants $\xL_2, \xL_4$ depend only on
the space dimension $N$ \cite{ggm03}.

These type of inequalities arise very naturally in the study of
singular differential operators. We would like to mention in
particular that  improved Hardy inequalities arise in the study of
singular solutions of the Gelfand problem \cite{bt, bv}, whereas
the improved Rellich in the biharmonic analogue of the Gelfand
problem \cite{ggm03}. It is worth  noting the work of Eilertsen
\cite{eil01} which is connected with the work of Maz'ya [M2] on
the Wiener test for higher order Elliptic equations. Related are
also the works of Yafaev \cite{yaf99} and Grillo \cite{gg}. For
some recent results concerning Hardy-Sobolev inequalities we refer
to \cite{a, acr, hn, ms, v}.

%
Our aim in this paper is to obtain sharp improved versions of
inequalities  such as (\ref{rel1}) and (\ref{1.7}), where
additional non-negative terms are present in the respective
right-hand sides. At the same time we obtain some new improved
Rellich inequalities which are new even at the level of plain
Rellich inequalities. The method we use was first introduced in
\cite{ft01} to obtain Hardy inequalities, here we extend it to
obtain higher order Rellich inequalities. Attached to the Rellich
inequality (\ref{rel1}), there is a similar Rellich inequality
that connects first to second order derivatives. That is, for $N
\geq 5,$ and for all $ u \in \RN$ we have
\be \label{rel2} \int_{\RN} (\xD u)^2 \, dx \geq \frac{N^2}{4}\,
\int_{\RN} \frac{|\ana u|^2}{|x|^2} \, dx. \ee
The constant $\frac{N^2}{4}$ is the best constant for
(\ref{rel2}). From this inequality and from (\ref{1.5}) we easily
arrive to a much stronger inequality than (\ref{rel1}). It was a
surprise for us that we have not trace inequality (\ref{rel2}) in
the literature.

From now on $\XO$ is a bounded domain containing the origin. In
$\XO$ inequalities (\ref{rel1}) and (\ref{rel2}) take the
following much stronger form.

\begin{theorem}\label{a.1}
{\bf (Improved Rellich-Sobolev inequality)} Let $N\geq 5$ and $D
\geq \sup_{x \in \XO} |x|$. There exists a positive constant c
such that for all $u \in H_0^2(\XO)$ there holds
\begin{eqnarray}
\ib&& \int_{\Omega}(\Delta u)^2dx \geq \biggl
(\frac{N(N-4)}{4}\biggr )^2\, \int_{\XO} \frac{u^2}{|x|^4} \, dx
+\, c\, \biggl (\int_{\Omega}|u|^{\frac{2N}{N-4}}\,
X_1^{\frac{2(N-2)}{N-4}}\, dx\biggr )^{\frac{N-4}{N}} \; .
\label{4.11} \\
\ia && \int_{\Omega}(\Delta u)^2dx\geq
\frac{N^2}{4}\,\int_{\Omega}\frac{|\nabla u|^2}{|x|^2} \, dx +c\,
\biggl (\int_{\Omega}|\nabla u|^{\frac{2N}{N-2}}\,
X_1^{\frac{2(N-1)}{N-2}}\, dx\biggr )^{\frac{N-2}{N}} \;
,\label{4.20}
\end{eqnarray}
\end{theorem}

Let us now give the following
\begin{defin}({\it Optimal Inequality})
Suppose that for some potential $V,$ we have for all $u \in \CO\,
,$ \be \label{1.11} \int_{\Omega}(\Delta u)^2dx\geq
\frac{N^2}{4}\,\int_{\Omega}\frac{|\nabla u|^2}{|x|^2} \, dx +\,
\int_{\Omega}|\nabla u|^2 \, V\, dx \;. \ee We say that inequality
(\ref{1.11}) is {\bf optimal}, when there is no potential $W
\gneqq 0$ to make the inequality \be \label{1.12}
\int_{\Omega}(\Delta u)^2dx\geq
\frac{N^2}{4}\,\int_{\Omega}\frac{|\nabla u|^2}{|x|^2} \, dx +\,
\int_{\Omega}|\nabla u|^2 \, V\, dx \;+\, \int_{\Omega}|\nabla
u|^2 \, W\, dx \;, \ee hold true for all $u \in \CO$.
\end{defin}

We then have
\begin{theorem}\label{a.2}
{\bf (Improved Rellich inequality I)} Let $N\geq 5$ and $D \geq
\sup_{x \in \XO} |x|$. \\
\ia \,Suppose the potential $V \gneqq 0$ is such that \be
\label{1.13} \int_{\Omega} V^{\frac{N}{2}} \, X_1^{1-N} \, dx \;\;
< \;\; +\infty. \ee Then there exists a positive constant b such
that for all $u \in \CO$, there holds \be \label{1.14}
\int_{\Omega}(\Delta u)^2dx\geq
\frac{N^2}{4}\,\int_{\Omega}\frac{|\nabla u|^2}{|x|^2} \, dx
\,+\,b\, \int_{\Omega}|\nabla u|^2 \, V\, dx \,. \ee If in
addition b is the best constant, then inequality (\ref{1.14}) is
an optimal inequality. \\
\ib \,Suppose the potential $W \gneqq 0$ is such that \be
\label{1.15} \int_{\Omega} W^{\frac{N}{4}} \, X_1^{1-\frac{N}{2}}
\, dx \;\; < \;\; +\infty. \ee Then there exists a positive
constant c such that for all $u \in \CO$, there holds \be
\label{1.16} \int_{\Omega}(\Delta u)^2dx\geq \biggl
(\frac{N(N-4)}{4}\biggr )^2\,\int_{\Omega}\frac{|u|^2}{|x|^4} \,
dx \,+\,c\, \int_{\Omega}|u|^2 \, V\, dx \,. \ee If in addition c
is the best constant, then inequality (\ref{1.16}) is an optimal
inequality.
\end{theorem}
The difficult part in the previous theorem is establishing that
inequalities are optimal and we will do that in section 4. We can
improve Rellich inequality differently and obtain
\begin{theorem} \label{t5.4}
{\bf (Improved Rellich inequality II)} Let $N\geq 5$ and $D \geq
\sup_{x \in \XO} |x|$. Then for all $u \in \CO,$ there holds
\be \label{5.8} \int_{\Omega} (\Delta u)^2\, dx \; \geq \; \biggl
(\frac{N(N-4)}{4} \biggr)^2\, \int_{\Omega} \frac{u^2}{|x|^4}\, dx
\; +\; \biggl ( 1 + \frac{N(N-4)}{8} \biggr ) \sum_{i=1}^{\infty}
\int_{\Omega} \frac{u^2}{|x|^4}\, X_{1}^{2} X_{2}^{2} \ldots
X_{i}^{2} \, dx. \ee
Moreover, for each $k=1,2,\ldots$,\ the constant\ $\biggl ( 1 +
\frac{N(N-4)}{8} \biggr )$\ is the best constant for the
corresponding k-Improved Rellich Inequality, that is
\begin{eqnarray} \label{ray1}
1 + \frac{N(N-4)}{8} =
\;\;\;\;\;\;\;\;\;\;\;\;\;\;\;\;\;\;\;\;\;\;\;\;\;\;\;\;\;\;\;
\;\;\;\;\;\;\;\;\;\;\;\;\;\;\;\;\;\;\;\;\;\;\;\;\;\;\;\;\;\;\;\;\;\;
\;\;\;\;\;\;\;\;\;\;\;\;\;\;\;\;\;\;\;\;\;\;\;\;\;\;\;\;\;\;\;\;\;\;
\;\;\;\;\;\;\;\;\;\;\;\;\;\;\;\;\;\nonumber \\
\inf_{u \in \CO} \frac{\int_{\Omega} (\Delta u)^2\, dx - \biggl (
\frac{N(N-4)}{4} \biggr )^2\, \int_{\Omega} \frac{u^2}{|x|^4}\, dx
-\biggl ( 1 + \frac{N(N-4)}{8} \biggr ) \sum_{i=1}^{k-1}
\int_{\xO} \frac{u^2}{|x|^4}
 X_1^2 X_2^2 \ldots X_i^2\;  dx}
{\int_{\xO}  \frac{u^2}{|x|^4} X_1^2 X_2^2
 \ldots X_k^2\;  dx}.
\end{eqnarray}
\end{theorem}
\begin{theorem} \label{t5.5}
{\bf (Improved Rellich inequality III)} Let $N\geq 5$ and $D \geq
\sup_{x \in \XO} |x|$. Then for all $u \in \CO,$ there holds 
\begin{equation} \label{5.9}
\int_{\Omega} (\Delta u)^2\, dx \geq  \frac{N^2}{4} \int_{\Omega}
\frac{|\nabla u|^2}{|x|^2}\, dx +\, \frac{1}{4}
\sum_{i=1}^{\infty} \int_{\Omega} \frac{|\nabla u|^2}{|x|^2}\,
X_{1}^{2} X_{2}^{2} \ldots X_{i}^{2} \, dx.
\end{equation}
Moreover, the constant\ $\frac{N^2}{4}$\ is the best and similarly
for each $k=1,2,\ldots$,\ the constant\ $\frac{1}{4}$\ is the best
constant for the corresponding k-Improved Rellich Inequality, that
is
\begin{equation}\label{ray2}
\frac{1}{4} = \inf_{u \in \CO} \frac{\int_{\Omega} (\Delta u)^2\,
dx - \frac{N^2}{4} \int_{\Omega} \frac{|\nabla u|^2}{|x|^2}\, dx
-\frac{1}{4} \sum_{i=1}^{k-1} \int_{\xO} \frac{|\nabla
u|^2}{|x|^2}
 X_1^2 X_2^2 \ldots X_i^2\;  dx}
{\int_{\xO}  \frac{|\nabla u|^2}{|x|^2} X_1^2 X_2^2
 \ldots X_k^2\;  dx}.
\end{equation}
\end{theorem}

Next we consider higher order Rellich inequalities.  When applying
Theorems \, \ref{t5.5}, \, \ref{t5.4} we reduce the order by one
or two.  In doing so weights enter in our inequalities. For this
reason we first consider second order Rellich inequalities with
weights. For\ $N\geq5$\ and\ $0 \leq m < \frac{N-4}{2}$,\ holds
that
\begin{equation} \label{rel1m}
\int_{\Omega} \frac{|\Delta u|^2}{|x|^{2m}}\, dx \geq \biggl (
\frac{(N+2m)(N-4-2m)}{4} \biggr )^2\, \int_{\Omega}
\frac{u^2}{|x|^{2m+4}}\, dx,
\end{equation}
while the corresponding improved inequalities can be stated as
\begin{theorem} \label{t8.4}
{\bf (Improved Rellich inequality IV)} Suppose $N\geq5$, $0 \leq m
< \frac{N-4}{2}$ and \ $D \geq \sup_{x \in \Omega} |x|$. Then for
all $ u \in \CO,$ there holds
\begin{eqnarray} \label{8.4}
\int_{\Omega} \frac{(\Delta u)^2}{|x|^{2m}}\, dx \, \geq \, \biggl
( \frac{(N+2m)(N-4-2m)}{4} \biggr )^2\, \int_{\Omega}
\frac{u^2}{|x|^{2m+4}}\, dx
\;\;\;\;\;\;\;\;\;\;\;\;\;\;\;\;\;\;\;\;\;\;\;\;\;\;\;\;\;\;
\;\;\;\;\;\;\;\;\;\;\;\;
\nonumber \\
\;\;\;\;\;\;\;\;+\, \biggl ( (1+m)^2 + \frac{(N+2m)(N-4-2m)}{8}
\biggr ) \sum_{i=1}^{\infty} \int_{\Omega}
\frac{u^2}{|x|^{2m+4}}\, X_{1}^{2} X_{2}^{2} \ldots X_{i}^{2} \,
dx.
\end{eqnarray}
Moreover $(\frac{(N+2m)(N-4-2m)}{4})^2$\ is the best constant.
Similarly for each $k=1,2,\ldots$,\ the constant\ $(1+m)^2 +
\frac{(N+2m)(N-4-2m)}{8}$\ is the best constant for the
corresponding k-Improved Hardy-Rellich Inequality, that is
\begin{eqnarray} \label{eqbestm1}
(1+m)^2 + \frac{(N+2m)(N-4-2m)}{8}
=\;\;\;\;\;\;\;\;\;\;\;\;\;\;\;\;\;\;\;\;\;\;
\;\;\;\;\;\;\;\;\;\;\;\;\;\;\;\;\;\;\;\;\;\;\;\;\;\;\;\;\;\;\;\;\;
\;\;\;\;\;\;\;\;\;\; \;\;\;\; \nonumber
\\
\inf_{u \in \CO} \frac{\int_{\Omega} \frac{|\Delta
u|^2}{|x|^{2m}}\, dx - \biggl ( \frac{(N+2m)(N-4-2m)}{4} \biggr
)^2 \int_{\Omega} \frac{u^2}{|x|^{2m+4}}\, dx - A \sum_{i=1}^{k-1}
\int_{\xO} \frac{u^2}{|x|^{2m+4}}
 X_1^2 \ldots X_i^2\; u^2 dx}
{\int_{\xO}  \frac{u^2}{|x|^{2m+4}} X_1^2
 \ldots X_k^2\;  dx},
\end{eqnarray}
where\ $A = (1+m)^2 + \frac{(N+2m)(N-4-2m)}{8}$.\
\end{theorem}

On the other hand the weighted Rellich inequality of the form
(\ref{rel2}) reads:
\begin{theorem} \label{a.6}
Suppose\ $N \geq 5$\ and\ $0 \leq m<\frac{N-4}{2}$.\ Then, for all
$u \in C_{0}^{\infty} (\Omega)$,\ there holds
\begin{equation} \label{rel2m}
\int_{\Omega} \frac{|\Delta u|^2}{|x|^{2m}}\, dx \geq a_{m,N}\,
\int_{\Omega} \frac{|\nabla u|^2}{|x|^{2m+2}}\, dx,
\end{equation}
where the best constant $a_{m,N}$\ is given by:
\begin{equation} \label{1.23}
a_{m,N} := \min_{k = 0,1,2,...} \frac{\biggl (
\frac{(N-4-2m)(N+2m)}{4}+k(N+k-2) \biggr )^2}{\biggl (
\frac{N-4-2m}{2} \biggr )^2 +k(N+k-2)}.
\end{equation}
In particular when\ $0 \leq m \leq \frac{-(N+4)+2\sqrt{N^2 -N
+1}}{6}$, we have
\[
a_{m,N} = \biggl ( \frac{N+2m}{2} \biggr )^2,
\]
whereas when\ $\frac{-(N+4)+2\sqrt{N^2 -N +1}}{6} < m <
\frac{N-4}{2}$, we have
\[
0 < a_{m,N} < \biggl ( \frac{N+2m}{2} \biggr )^2.
\]
\end{theorem}
In Theorem \ref{p9.1} we have a full description of how the
constant $a_{m,N}$ behaves. Our next result is

\begin{theorem} \label{t9.2}
{\bf (Improved Rellich inequality V)} Let\ $D \geq \sup_{x \in
\Omega} |x|$\ and\ $0\leq m \leq \frac{-(N+4)+2\sqrt{N^2 -N
+1}}{6}$.\ Then for all $u \in \CO$, there holds
\begin{eqnarray} \label{9.40}
\int_{\Omega} \frac{|\Delta u|^2}{|x|^{2m}}\, dx - \biggl (
\frac{N+2m}{2} \biggr )^2\, \int_{\Omega} \frac{|\nabla
u|^2}{|x|^{2m+2}}\, dx \geq \frac{1}{4}   \sum_{i=1}^{\infty}
\int_{\Omega} \frac{|\nabla u|^2}{|x|^{2m+2}}\, X_{1}^{2}
X_{2}^{2} \ldots X_{i}^{2} \, dx.
\end{eqnarray}
Moreover for each $k=1,2,\ldots$,\ the constant\ $\frac{1}{4}$\ is
the best constant for the corresponding k-Improved Hardy-Rellich
Inequality, that is
\begin{equation} \label{e7.7}
\frac{1}{4} = \inf_{u \in H_0^2 (\xO)} \frac{\int_{\Omega}
\frac{|\Delta u|^2}{|x|^{2m}}\, dx - \biggl ( \frac{N+2m}{2}
\biggr )^2 \int_{\Omega} \frac{|\nabla u|^2}{|x|^{2m+2}}\, dx
-\frac{1}{4} \sum_{i=1}^{k-1} \int_{\xO} \frac{|\nabla
u|^2}{|x|^{2m+2}}
 X_1^2 X_2^2 \ldots X_i^2\; dx}
{\int_{\xO}  \frac{|\nabla u|^2 }{|x|^{2m+2}} X_1^2 X_2^2
 \ldots X_k^2\; dx}.
\end{equation}
\end{theorem}

In order to state our improved higher order Rellich inequality we
set
\begin{equation} \label{1.26}
\sigma(m,N) = \biggl ( \frac{(N+2m)(N-4-2m)}{4} \biggr )^2 \ee \be
\label{1.27}\bar{\sigma}(m,N) = (1+m)^2 +
\frac{(N+2m)(N-4-2m)}{8}.
\end{equation}
We then have
\begin{theorem} \label{a.8}
{\bf (Improved Higher Order Rellich Inequalities I)} Suppose $m
\in \mathbb{N},\,l=0, \cdots\, ,m-1,\; 4m < N$ and $D \geq \sup_{x
\in \Omega} |x|$. Then for all $u \in C_{0}^{\infty}(\Omega)$
there holds
\begin{eqnarray}
\ia \; \int_{\Omega} (\Delta^m u)^2\, dx &\geq&
\prod_{k=0}^{l}\;\biggl ( \frac{(N+4k)(N-4-4k)}{4} \biggr )^2\;
\int_{\Omega} \frac{(\Delta^{m-l-1} u)^2}{|x|^{4l+4}}\, dx  \nonumber \\
&+&\sum_{k=1}^{l}\bar{\sigma}(2k,N)\; \prod_{j=0}^{k-1}\sigma
(2j,N) \sum_{i=1}^{\infty} \int_{\Omega} \frac{(\Delta^{m-k-1}
u)^2}{|x|^{4k+4}} X_{1}^{2}
\ldots X_{i}^{2}dx\   \nonumber \\
&+&\biggl (1+\, \frac{N(N-4)}{8} \biggr )\;\sum_{i=1}^{\infty}\;
\int_{\Omega} \frac{(\Delta^{m-1} u)^2}{|x|^{4}}\, X_{1}^{2}
\ldots
X_{i}^{2}\,dx \; , \label{1.28} \\
\nonumber \\
 \ib \int_{\Omega} |\nabla \Delta^m u|^2\,
dx &\geq& \; \biggl ( \frac{N-2}{2} \biggr )^2\;
\prod_{k=0}^{l-1}\;\biggl ( \frac{(N+2+4k)(N-6-4k)}{4} \biggr
)^2\, \int_{\Omega} \frac{(\Delta^{m-l}
u)^2}{|x|^{4l+2}}\, dx  \nonumber \\
&+&  \biggl ( \frac{N-2}{2} \biggr )^2 \;
\sum_{k=2}^{l}\bar{\sigma}(2k-1,N)\; \prod_{j=0}^{k-2}\;\sigma
(2j+1,N)\, \sum_{i=1}^{\infty}\; \int_{\Omega} \frac{(\Delta^{m-k}
u)^2}{|x|^{4k+2}}\, X_{1}^{2} \ldots
X_{i}^{2}\,dx \nonumber \\
&+& \biggl ( \frac{N-2}{2} \biggr )^2 \;\biggl (4+
\frac{(N+2)(N-6)}{8} \biggr )
 \, \sum_{i=1}^{\infty}\;
\int_{\Omega} \frac{(\Delta^{m-1} u)^2}{|x|^{6}}\, X_{1}^{2}
\ldots
X_{i}^{2}\,dx \nonumber \\
&+&  \frac{1}{4}\;  \sum_{i=1}^{\infty}\; \int_{\Omega}
\frac{(\Delta^{m} u)^2}{|x|^{2}}\, X_{1}^{2} \ldots X_{i}^{2}\,dx
\; . \label{1.29}
\end{eqnarray}
\end{theorem}
\begin{theorem} \label{tp}
{\bf (Improved Higher Order Rellich Inequality II)} Suppose $m, \,
l \, \in \mathbb{N}$,\ $1 \, \leq \, l \, \leq \, \frac{-N+8+2
\sqrt{N^2-N+1}}{12}$,\ $4m < N$ and $D \geq \sup_{x \in \Omega}
|x|$. Then for all $u \in C_{0}^{\infty}(\Omega)$ there holds
\begin{eqnarray}
\int_{\Omega} (\Delta^m u)^2\, dx &\geq& \prod_{k=0}^{l-1}\;\biggl
( \frac{(N+4k)(N-4-4k)}{4} \biggr )^2\; \int_{\Omega}
\frac{(\Delta^{m-l}
u)^2}{|x|^{4l}}\, dx  \nonumber \\
&+& \frac{4}{(N-4)^2} \sum_{k=1}^{l}\; \prod_{j=0}^{k-1}\;\biggl (
\frac{(N-4+4j)(N-4j)}{4} \biggr )^2\; \sum_{i=1}^{\infty}\;
\int_{\Omega} \frac{(\nabla \,\Delta^{m-k} u)^2}{|x|^{4k-2}}\,
X_{1}^{2} \ldots X_{i}^{2}\,dx\;\nonumber \\
&+& \frac{1}{N^2} \sum_{k=1}^{l}\; \prod_{j=1}^{k}\;\biggl (
\frac{(N+4j)(N-4j)}{4} \biggr )^2\; \sum_{i=1}^{\infty}\;
\int_{\Omega} \frac{(\Delta^{m-k} u)^2}{|x|^{4k}}\, X_{1}^{2}
\ldots X_{i}^{2}\,dx\, . \label{1.30}
\end{eqnarray}
\end{theorem}

The paper is divided in two parts. In the first part we deal with
the biharmonic operator, while in the second part we deal with the
polyharmonic operator. More precisely, in Section 2 we prove some
identities and inequalities to be used widely in the sequel; the
main tool for this is decomposition into spherical harmonics. In
Section 3 we prove Theorems \ref{a.1}, \ref{t5.4} and \ref{t5.5},
while in Section 4 we prove that the constants appearing in
certain inequalities are the best and complete the proof of
Theorems \ref{t5.4} and \ref{t5.5}. In Section 5 we state
necessary conditions for the improvement or not of inequalities
(\ref{rel1}), (\ref{rel2}). In the two Sections of Part II we
actually prove Theorems \ref{t8.4} to \ref{tp}. \vspace{0.2cm}

\emph{Notation}: Sometimes, for the sake of the representation we
use the following quantities
\begin{eqnarray*}
I_{\Omega}[u] &:=& \int_{\Omega} |\Delta u|^2\, dx - \biggl (
\frac{N(N-4)}{4} \biggr )^2\, \int_{\Omega} \frac{u^2}{|x|^4}\,
dx, \\
J_{\Omega}[v] &:=& \int_{\Omega} |x|^{-(N-4)} |\Delta v|^2\, dx -
N(N-4) \int_{\Omega} |x|^{-N} (x \cdot \nabla v)^2\, dx \\ && +
\frac{N(N-4)}{2} \int_{\Omega} |x|^{-(N-2)} |\nabla v|^2\, dx, \\
\mathbb{I}_{\Omega}[u] &:=& \int_{\Omega} |\Delta u|^2\, dx -
\frac{N^2}{4} \int_{\Omega} \frac{|\nabla u|^2}{|x|^2}\, dx, \\
\mathbb{J}_{\Omega} [v] &:=& \int_{\Omega} |x|^{-(N-4)} |\Delta
v|^2\, dx - N(N-4) \int_{\Omega} |x|^{-N} (x \cdot \nabla v)^2\,
dx \nonumber
\\ && + \frac{N(N-8)}{4} \int_{\Omega} |x|^{-(N-2)} |\nabla v|^2\,
dx,
\end{eqnarray*}
related to (\ref{rel1}) and (\ref{rel2}). \vspace{0.5cm}
%
%
%
\begin{center}
{\bf PART I. THE BIHARMONIC OPERATOR}
\end{center}
%
%
\section{Preliminaries}
\setcounter{equation}{0}
In this section we establish some abstract relations to be used in
the sequel. In the first part we prove some useful identities
while, in the second part we apply spherical harmonic
decomposition in order to prove certain inequalities. Throughout
this section\ $\Omega$\ is an arbitrary domain (bounded or
unbounded).
\subsection{Preliminaries Identities}
\setcounter{equation}{0}
\begin{lemma} \label{l2.1a}
Let\ $N \geq 3$,\ $a < N-2$\ and\ $B \in C^2 [0,+\infty)$.\ Then,
for any\ $u \in C_{0}^{\infty} (\Omega)$,\ we have the identity:
\[
\int_{\Omega} \frac{B(r)}{r^a} |\nabla u|^2\, dx = - \int_{\Omega}
\frac{B(r)}{r^a} u\, \Delta u\, dx + \frac{1}{2} \int_{\Omega}
\Delta \biggl ( \frac{B(r)}{r^a} \biggr )\, u^2\, dx.
\]
\end{lemma}
\emph{Proof}\, Observe that for\ $a < N-2$\ we have that\ $\Delta
\biggl ( \frac{B(r)}{r^a} \biggr )\, u^2 \in L^1 (\Omega)$\ and\
$\nabla \biggl ( \frac{B(r)}{r^a} \biggr )\, \nabla u^2 \in L^1
(\Omega)$.\ In virtue of the identity
\begin{equation} \label{2.1aa}
|\nabla w|^2 = \frac{1}{2} \Delta w^2 - w \Delta w,
\end{equation}
it suffices to prove that
\begin{equation} \label{2.1a}
\int_{\Omega} \frac{B(r)}{r^a}\, \Delta u^2\, dx = \int_{\Omega}
\Delta \biggl ( \frac{B(r)}{r^a} \biggr )\, u^2\, dx.
\end{equation}
If we write
\[
\int_{\Omega} \frac{B(r)}{r^a}\, \Delta u^2\, dx = \int_{\Omega
\backslash B_{\varepsilon}} \frac{B(r)}{r^a}\, \Delta u^2\, dx +
\int_{B_{\varepsilon}} \frac{B(r)}{r^a}\, \Delta u^2\, dx,
\]
and using the limits
\[
\bigg | \int_{\partial B_{\varepsilon}} \frac{\partial}{\partial
\nu} \biggl ( \frac{B(r)}{r^a} \biggr ) u^2\, ds \bigg | \leq c\;
\varepsilon^{N-2-a} \to 0\;\;\; \mbox{and}\;\;\; \bigg |
\int_{\partial B_{\varepsilon}} \frac{B(r)}{r^a}
\frac{\partial}{\partial \nu} u^2\, ds \bigg | \leq c\;
\varepsilon^{N-1-a} \to 0,
\]
as\ $\varepsilon \to 0,$\ we obtain that (\ref{2.1a}) is true.
Thus, the proof is completed.\ $\blacksquare$
\begin{lemma} \label{lemma2.1}
Let\ $N>4$\ and\ $0 < a \leq \frac{N-4}{2}$.\ For any\ $u \in
C_{0}^{\infty} (\Omega)$\ we set\ $v = |x|^{a} u$.\ Then, the
following equality holds.
\begin{eqnarray} \label{2.1}
\int_{\Omega} |\Delta u|^2\, dx = \int_{\Omega} |x|^{-2a} |\Delta
v|^2\, dx - 4a(a+2) \int_{\Omega} |x|^{-2a-4} (x \cdot \nabla
v)^2\, dx \nonumber \\ + 2a(a+2) \int_{\Omega} |x|^{-2a-2} |\nabla
v|^2\, dx \nonumber \\ + a(a+2)(-N+a+2)(-N+a+4) \int_{\Omega}
|x|^{-2a-4} v^2\, dx,
\end{eqnarray}
\end{lemma}
\emph{Proof}\, We have that
\begin{eqnarray}\label{2.2}
\int_{\Omega} | \Delta u |^2\, dx &=& \int_{\Omega} |x|^{-2a}
|\Delta v |^2\, dx + 4a^2 \int_{\Omega} |x|^{-2a-4} (x \cdot
\nabla v)^2\, dx \nonumber \\ && +  \biggl ( -aN+a(a+2) \biggr )^2
\int_{\Omega} |x|^{-2a-4} v^2\, dx + I_1 + I_2 + I_3,
\end{eqnarray}
where
\begin{eqnarray*}
I_1 &:=& 2 \int_{\Omega} |x|^{-a} \Delta v\, v \Delta |x|^{-a}\,
dx \\
I_2 &:=& 4 \int_{\Omega} |x|^{-a} \Delta v\, \nabla v \cdot \nabla
|x|^{-a}\, dx \\
I_3 &:=& 4 \int_{\Omega} \nabla |x|^{-a} \cdot \nabla v\, v \Delta
|x|^{-a}\, dx.
\end{eqnarray*}
Following the same procedure as in the proof of Lemma \ref{l2.1a}
(having in mind also that\ $|v| \leq |x|^{a} ||u||_{\infty}$) we
obtain that
\begin{eqnarray*}
I_1 &=& a(-N+a+2)(2a+2)(-N+2a+4) \int_{\Omega} |x|^{-2a-4} v^2\,
dx \nonumber \\ &&- 2 \biggl ( -a N +a(a+2) \biggr ) \int_{\Omega}
|x|^{-2a-2} | \nabla v|^2\, dx,  \\
I_2 &=& 4a(-2a-2) \int_{\Omega} |x|^{-2a-4} (x \cdot \nabla v)^2\,
dx + 2a (-N+2a+4) \int_{\Omega} |x|^{-2a-2} |\nabla v|^2\, dx,  \\
I_3 &=& 2a^2 (-N+a+2)(-N+2a+4) \int_{\Omega} |x|^{-2a-4} v^2\, dx.
\end{eqnarray*}
Then, from (\ref{2.2}) we conclude (\ref{2.1}) and the proof is
completed.\ $\blacksquare$ \vspace{0.2cm} \\
Using the previous lemma we may easily obtain the following
result, concerning the relation between\ $I, \mathbb{I}, J,
\mathbb{J}$.\
\begin{lemma} \label{l2.2}
Let\ $N\geq5$,\ $u \in C_{0}^{\infty} (\Omega)$\ and\ $v =
|x|^{(N-4)/2} u$.\ We have that:
\begin{eqnarray*}
i)\;\;\; \int_{\Omega} \frac{|\nabla u|^2}{|x|^2}\, dx =
\int_{\Omega} |x|^{-(N-2)} |\nabla v|^2\, dx + \biggl (
\frac{N-4}{2} \biggr )^2 \int_{\Omega} |x|^{-N} |v|^2\,
dx,\;\;\;\;\;\;\;\;\;\;\;\;\;\;\;\;\;\;\;\;
\end{eqnarray*}
\begin{eqnarray*}
ii)\;\;\; \int_{\Omega} |\Delta u|^2\, dx - \biggl (
\frac{N(N-4)}{4} \biggr )^2\, \int_{\Omega} \frac{u^2}{|x|^4}\, dx
= \int_{\Omega} |x|^{-(N-4)} |\Delta v|^2\, dx -
\;\;\;\;\;\;\;\;\;\;\;\;\;\;\;\;\;\;\;\;\;\;\;\;\;\;\;\;\;\;\;\\
- N(N-4) \int_{\Omega} |x|^{-N} (x \cdot \nabla v)^2\, dx  +
\frac{N(N-4)}{2} \int_{\Omega} |x|^{-(N-2)} |\nabla v|^2\, dx,
\end{eqnarray*}
\begin{eqnarray*}
iii)\;\;\; \int_{\Omega} |\Delta u|^2\, dx - \frac{N^2}{4}
\int_{\Omega} \frac{|\nabla u|^2}{|x|^2}\, dx = \int_{\Omega}
|x|^{-(N-4)} |\Delta v|^2\, dx -
\;\;\;\;\;\;\;\;\;\;\;\;\;\;\;\;\;\;\;\;\;\;\;\;\;\;\;\;\;\;\;\;\;\;\;
\;\;\;\;\;\;\;\;\;\;\\
- N(N-4) \int_{\Omega} |x|^{-N} (x \cdot \nabla v)^2\, dx +
\frac{N(N-8)}{4} \int_{\Omega} |x|^{-(N-2)} |\nabla v|^2\, dx,
\end{eqnarray*}
\end{lemma}
%

%
%
%
\subsection{Preliminaries Inequalities}
The decomposition of\ $u$\ and\ $v$\ into spherical harmonics will
be one of the main tools in our investigation. Let\ $u \in
C_{0}^{\infty} (\Omega)$.\ If we extend\ $u$\ as zero outside\
$\Omega$,\ we may consider that\ $u \in C_{0}^{\infty}
(\mathbb{R}^N)$.\ Decomposing\ $u$\ into spherical harmonics we
get
\[
u = \sum_{k=0}^{\infty} u_k := \sum_{k=0}^{\infty} f_k (r) \phi_k
(\sigma),
\]
where\ $\phi_k (\sigma)$\ are the orthonormal eigenfunctions of
the Laplace-Beltrami operator with corresponding eigenvalues\ $c_k
= k(N+k-2)$,\ $k \geq 0$.\ The functions\ $f_k$\ belong in\
$C^{\infty}_{0} (\Omega)$, satisfying\ $f_k (r) = O (r^k)$\ and\
$f_k' (r) = O (r^{k-1})$,\ as\ $r \downarrow 0$.\ In particular,\
$\phi_0 (\sigma) =1$\ and\ $u_0 (r) = \frac{1}{|\partial B_r|}
\int_{\partial B_r} u\, ds$,\ for any\ $r>0$.\ Then, for any\ $k
\in \mathbb{N}$,\ we have that
\[
\Delta u_k = \biggl ( \Delta f_k (r) - \frac{c_k f_k(r)}{r^2}
\biggr ) \phi_k (\sigma)
\]
so
\begin{equation} \label{2.3}
\int_{\mathbb{R}^N} |\Delta u_k|^2\, dx = \int_{\mathbb{R}^N}
\biggl ( \Delta f_k (r) - \frac{c_k f_k(r)}{r^2} \biggr )^2\, dx.
\end{equation}
In addition,
\begin{equation} \label{2.4}
\int_{\mathbb{R}^N} |\nabla u_k|^2\, dx = \int_{\mathbb{R}^N}
\biggl ( |\nabla f_k (r)|^2 + c_k \frac{f^2_k(r)}{r^2} \biggr )\,
dx.
\end{equation}

Next, we assume the function\ $v \in C_0^{\infty} (\mathbb{R}^N
\backslash \{0\})$,\ such that\ $v = |x|^{\frac{N-4}{2}} u$.\ From
the definitions of\ $u$\ and\ $v$,\ we may write that
\[
u = \sum_{k=0}^{\infty} u_k = \sum_{k=0}^{\infty}
r^{\frac{-N+4}{2}+k} g_k (r) \phi_k (\sigma),\;\;\;
v = \sum_{k=0}^{\infty} v_k = \sum_{k=0}^{\infty} r^{k} g_k (r)
\phi_k (\sigma),
\]
where\ $f_k = r^{-\frac{N-4}{2}+k} g_k$,\ with\ $g_k \sim 0$\ and\
$r\, g'_k \sim 0$\ at the origin. More precisely, we may prove
that the following identities hold, for any\ $k \in \mathbb{N}$.\
\begin{eqnarray} \label{3.Du}
\int_{\mathbb{R}^N} |\Delta u_k|^2 dx = \int_{\mathbb{R}^N}
r^{2k-N+4} |\nabla g^{'}_k|^2 dx + \biggl ( \frac{N(N-4)}{2} + 2k
(N-3) +3 \biggr ) \int_{\mathbb{R}^N} r^{2k-N+2} |\nabla g_k|^2 dx
\nonumber \\
+ \biggl [ \biggl ( \frac{N(N-4)}{4} \biggr )^2 + \frac{N(N-4)}{2}
(c_k +k^2) \biggr ] \int_{\mathbb{R}^N} r^{2k-N} (g_k)^2\, dx
\end{eqnarray}
\begin{eqnarray} \label{3.Gu}
\int_{\mathbb{R}^N} |x|^{-2} |\nabla u_k|^2\, dx =
\int_{\mathbb{R}^N} r^{2k-N+2} |\nabla g_k|^2\, dx + \biggl [
\biggl ( \frac{N-4}{2} \biggr )^2 + k(N-2) \biggr ]
\int_{\mathbb{R}^N} r^{2k-N} (g_k)^2\, dx,
\end{eqnarray}
\begin{eqnarray} \label{3.I}
I[u_k] = \int_{\mathbb{R}^N} r^{2k-N+4} |\nabla g^{'}_k|^2\, dx +
\biggl ( \frac{N(N-4)}{2} + 2k(N-3) +3 \biggr )
\int_{\mathbb{R}^N} r^{2k-N+2} |\nabla g_k|^2\, dx\;\;\;\;\;\;\;
\;\;\; \nonumber \\
 + \biggl [ \frac{N(N-4)}{2} (c_k +k^2) \biggr ] \int_{\mathbb{R}^N} r^{2k-N}
(g_k)^2\, dx,\;\;\;\;\;\;\;\;
\end{eqnarray}
\begin{eqnarray} \label{3.II}
\mathbb{I}[u_k] = \int_{\mathbb{R}^N} r^{2k-N+4} |\nabla
g^{'}_k|^2\, dx + \biggl ( (2k+N-1)(N-3) - \frac{N(3N-8)}{4}
\biggr )
\int_{\mathbb{R}^N} r^{2k-N+2} |\nabla g_k|^2\, dx \nonumber \\
 + \biggl [ \frac{N(3N-8)}{4} k^2 + \frac{N(N-8)}{4} c_k \biggr ]
\int_{\mathbb{R}^N} r^{2k-N} (g_k)^2\, dx,\;\;\;\;\;\;\;\;
\end{eqnarray}
\begin{eqnarray} \label{2.8}
\int_{\mathbb{R}^N} r^{-(N-4)} |\Delta v_k|^2\, dx =
\int_{\mathbb{R}^N} r^{2k-N+4} |\nabla g^{'}_k|^2\, dx
 + (2k+N-1)(N-3) \int_{\mathbb{R}^N} r^{2k-N+2} |\nabla
g_k|^2\, dx,
\end{eqnarray}
\begin{eqnarray} \label{2.9}
\int_{\mathbb{R}^N} r^{-(N-2)} |\nabla v_k|^2\, dx =
\int_{\mathbb{R}^N} r^{2k-N+2} |\nabla g_k|^2\, dx + k (N-2)
\int_{\mathbb{R}^N} r^{2k-N} (g_k)^2\, dx,\;\;\;\;\;\;\;\;\;
\;\;\;\;\;\;\;\;\;\;\;
\end{eqnarray}
\begin{eqnarray} \label{2.10}
\int_{\mathbb{R}^N} r^{-N} (x \cdot \nabla v_k)^2\, dx =
\int_{\mathbb{R}^N} r^{2k-N+2} |\nabla g_k|^2\, dx - k^2
\int_{\mathbb{R}^N} r^{2k-N} (g_k)^2\, dx.\;\;\;\;\;\;\;\;\;
\;\;\;\;\;\;\;\;\;\;\;\;\;\;\;\;\;\;\;\;\;\;
\end{eqnarray}
\vspace{0.3cm} \\
Let\ $k =1,2,...$\ and\ $V(|x|) \in C^{1} ([0,+\infty))$.\ The
following relation
\[
\int_{\mathbb{R}^N} V(|x|)\, |x|^{-2} |\grad u_k|^2 dx =
\int_{\mathbb{R}^N} V(|x|)\, |x|^{-2} |\grad f_k|^2 dx + c_k
\int_{\mathbb{R}^N} V(|x|)\, |x|^{-4} f_k^2 dx.
\]
implies that
\begin{eqnarray} \label{2.12}
\int_{\mathbb{R}^N} V(|x|)\, |x|^{-2} |\grad u_k|^2\, dx =
\int_{\mathbb{R}^N} r^{2k+2-N} V(|x|)\, |\nabla g_k|^2\, dx +
\;\;\;\;\;\;\;\;\;
\;\;\;\;\;\;\;\;\;\;\;\;\;\;\;\;\;\;\;\;\;\;\;\;\;\;\;\;\;\;\;\;
\;\;\;\;\;\;\; \nonumber \\
+ \biggl [ \biggl ( \frac{N-4}{2} \biggr )^2 +k(N-2) \biggr ]
\int_{\mathbb{R}^N}r^{2k-N} V(|x|)\, g_k^2\, dx + \biggl (
\frac{N-4}{2} -k \biggr ) \int_{\mathbb{R}^N} r^{2k+1-N} V'(|x|)\,
g_k^2\, dx.
\end{eqnarray}
Also, as an immediate consequence of the Hardy inequality, we have
the following relations.
\begin{eqnarray} \label{2.17}
\int_{0}^{\infty} r^{2k+3} (g_k'')^2\, dr &\geq& (k+1)^2
\int_{0}^{\infty} r^{2k+1} (g_k')^2\, dr, \\
\label{2.18} \int_{0}^{\infty} r^{2k+1} (g_k')^2\, dr &\geq& k^2
\int_{0}^{\infty} r^{2k-1} g_k^2\, dr.
\end{eqnarray}
Observe that in the case of a bounded domain\ $\Omega$,\ all the
obtained equalities remain true if we assume\ $B_D$,\ with\ $D =
sup_{x \in \Omega} |x|$,\ instead of\ $\mathbb{R}^N$.
\vspace{0.2cm}

In the remaining part of this section, using the decomposition
into spherical harmonics, we establish certain inequalities
concerning\ $I[u]$\ and\ $\mathbb{I}[u]$.\

\begin{theorem}
 \label{p3.1}
Let\ $N\geq 5$,\ $u \in C_{0}^{\infty} (\Omega)$\ and\ $v =
|x|^{(N-4)/2} u$.\ Then
\begin{equation} \label{3.1}
i)\;\;\; \int_{\Omega} |\Delta u|^2\, dx - \biggl (
\frac{N(N-4)}{4} \biggr )^2\, \int_{\Omega} \frac{u^2}{|x|^4}\, dx
\geq \biggl ( 4 + \frac{N(N-4)}{2} \biggr ) \int_{\Omega}
|x|^{-(N-2)} |\nabla v|^2\, dx.
\end{equation}
\begin{equation} \label{3.7}
ii)\;\;\; \int_{\Omega} |\Delta u|^2\, dx - \frac{N^2}{4}
\int_{\Omega} \frac{|\nabla u|^2}{|x|^2}\, dx \geq \biggl (
\frac{N-4}{2} \biggr )^2 \int_{\Omega} |x|^{-(N-2)} |\nabla v|^2\,
dx.
\end{equation}
\end{theorem}
\emph{Proof}\, i)\ It suffices to prove, by using (\ref{3.I}),
(\ref{2.9}) and (\ref{2.17}), that the following inequality
\begin{eqnarray}
\biggl [ (k+N-2)^2 - \frac{N(N-4)}{2} - \biggl ( 4 +
\frac{N(N-4)}{2} \biggr ) \biggr ] \int_{0}^{\infty} r^{2k+1}
(g'_k)^2\, dr \geq \nonumber \;\;\;\;\;\;\;\;\;\;\;\;
\;\;\;\;\;\;\;\;\;\;\;\; \\
\label{3.2}
\biggl [k (N-2) \biggl ( 4 + \frac{N(N-4)}{2} \biggr ) -
\frac{N(N-4)}{2} (k^2 + c_k) \biggr ] \int_{0}^{\infty} r^{2k-1}
(g_k)^2\, dr ,
\end{eqnarray}
holds for any\ $k=1,2,...$.\ or equivalently
\[
\biggl ( k+2N-4 \biggr ) \int_{0}^{\infty} r^{2k+1} (g'_k)^2\, dr
\geq \biggl [4 (N-2) -kN(N-4) \biggr ] \int_{0}^{\infty} r^{2k-1}
(g_k)^2\, dr ,
\]
which is true since
\[
k^2 \geq \frac{4 (N-2) -kN(N-4)}{k+2N-4},
\]
for\ $k=1,2,...$\ and\ $N\geq 5$.\ \vspace{0.2cm} \\
ii)\ From Lemma \ref{l2.2} we deduce that
\[
\mathbb{I}_{\Omega}[u] = I_{\Omega}[u] - \frac{N^2}{4}\,
\int_{\Omega} |x|^{-(N-2)} |\nabla v|^2\, dx.
\]
Then, the result follows from (\ref{3.1}).\ $\blacksquare$
\begin{lemma} \label{l3.2}
Let\ $N\geq 5$,\ $u \in C_{0}^{\infty} (\Omega)$\ and\ $v =
|x|^{(N-4)/2} u$.\ Then, the following inequalities hold.
\begin{eqnarray} \label{3.4}
\int_{\Omega} |x|^{4-N} |\Delta v|^2\, dx \geq N(N-4)
\int_{\Omega} |x|^{-N} (x \cdot \nabla v)^2\, dx + 4 \int_{\Omega}
|x|^{2-N} |\nabla v|^2\, dx
\end{eqnarray}
\begin{eqnarray} \label{3.3}
\int_{\Omega} |x|^{4-N} |\Delta v|^2\, dx \geq 2 (N-2)^2 \biggl (
\int_{\Omega} |x|^{-N} (x \cdot \nabla v)^2\, dx - \frac{1}{2}
\int_{\Omega} |x|^{2-N} |\nabla v|^2\, dx \biggr).
\end{eqnarray}
\end{lemma}
\emph{Proof}\, Inequality (\ref{3.4}) follows Theorem \ref{p3.1},
while (\ref{3.3}) follows from (\ref{3.4}) and the following
inequality
\begin{eqnarray} \label{3.5}
\int_{\Omega} |x|^{-N} |x \cdot \nabla v|^2\, dx - \frac{1}{2}
\int_{\Omega} |x|^{2-N} |\nabla v|^2\, dx \leq
\;\;\;\;\;\;\;\;\;\;\;\;\;\;\;\;\;\;\;\;\;\;\;\;\;\;\;\;\;\;\;
\;\;\;\;\;\;\;\;\;\;\;\;\;\;\;\;\;\;\;\;\;\;\;\;\;\;\;\;\;\;\;
\nonumber
\\ \frac{1}{2(N-2)^2} \biggl [ N(N-4) \int_{\Omega}
|x|^{-N} |x \cdot \nabla v|^2\, dx + 4 \int_{\Omega} |x|^{2-N}
|\nabla v|^2\, dx \biggr ].\;\;\; \blacksquare \vspace{0.1cm}
\end{eqnarray}
An immediate consequence of the inequality (\ref{3.3}) is the
following result.
%

\begin{coro} \label{p3.3}
Let\ $N\geq 5$,\ $u \in C_{0}^{\infty} (\Omega)$\ and\ $v =
|x|^{(N-4)/2} u$.\ Then
\begin{equation} \label{3.6}
\int_{\Omega} |\Delta u|^2\, dx - \biggl ( \frac{N(N-4)}{4} \biggr
)^2\, \int_{\Omega} \frac{u^2}{|x|^4}\, dx \geq \biggl (
\frac{1}{2} + \frac{2}{(N-2)^2} \biggr ) \int_{\Omega}
|x|^{-(N-4)} |\Delta v|^2\, dx.
\end{equation}
\end{coro}
\setcounter{theorem}{6}
\begin{theorem} \label{p3.5}
Let\ $N\geq 5$,\ $u \in C_{0}^{\infty} (\Omega)$\ and\ $v =
|x|^{(N-4)/2} u$.\ Then
\begin{equation} \label{3.8}
\int_{\Omega} |\Delta u|^2\, dx - \frac{N^2}{4} \int_{\Omega}
\frac{|\nabla u|^2}{|x|^2}\, dx \geq \biggl ( \frac{N-4}{2(N-2)}
\biggr )^2 \int_{\Omega} |x|^{-(N-4)} |\Delta v|^2\, dx.
\end{equation}
\end{theorem}
\emph{Proof}\, Using the identities (\ref{3.II}) and (\ref{2.8})
we have that (\ref{3.8}) holds if the following inequality
\begin{eqnarray*}
A := \int_{0}^{\infty} r^{2k+3} (g_k'')^2\, dr + (2kN-6k-1)
\int_{0}^{\infty} r^{2k+1} (g_k')^2\,
dr\;\;\;\;\;\;\;\;\;\;\;\;\;\;\; \nonumber \\ + (N-2)^2 (k^2 +
\frac{N-8}{3N-8} c_k) \int_{0}^{\infty} r^{2k-1} g_k^2\, dr \geq
0,
\end{eqnarray*}
is true for any\ $k \in \mathbb{N}$.\ Taking now into account
(\ref{2.17}) and (\ref{2.18}) we deduce that for any\ $k \in
\mathbb{N}$\ holds that
\[
A \geq A(k) \int_{0}^{\infty} r^{2k-1} g_k^2\, dr,
\]
where
\[
A(k) = k^2(k^2+2kN-4k)+k(N-2)^2(k + \frac{N-8}{3N-8}(k+N-2)).
\]
It is clear that\ $A(k)$\ is an increasing function for positive\
$k$,\ with\ $A(0)=0$.\ Thus, $A\geq0$,\ for any\ $k \in
\mathbb{N}$\ and the proof is completed.\ $\blacksquare$
%

%
%
\section{Hardy-Sobolev and Improved Inequalities}
\setcounter{equation}{0}
In this section we prove certain Hardy-Sobolev-type inequalities
and we establish some improved Hardy inequalities. Throughout this
section we assume that\ $N \geq 5$,\ $\Omega$\ is a bounded domain
and\ $D = sup_{x \in \Omega} |x|$.\ We extend any\ $u \in
C_{0}^{\infty} (\Omega)$\ as zero outside\ $\Omega$\ so we
consider that\ $u \in C_{0}^{\infty} (\mathbb{R}^N)$.\ We then
define\ $u_0 (r) := \frac{1}{|\partial B_r|} \int_{\partial B_r}
u\, ds$,\ for any\ $r>0$.\ It is clear that\ $u_0 \in
C_{0}^{\infty} [0,D)$.\
\begin{theorem} \label{p4.1}
Let\ $\Omega$\ be a bounded domain,\ $D = sup_{x \in \Omega} |x|$\
and\ $u \in C_{0}^{\infty} (\Omega)$.\ Then,
\begin{equation} \label{4.1}
I_{\Omega}[u] \geq I_{B_D}[u_0] + \frac{8(N-1)(N^2 -2N -2)}{(N^2
-4)^2} \int_{B_D} |\Delta (u-u_0)|^2\, dx.
\end{equation}
\end{theorem}
\emph{Proof}\, Observe that\ $I_{\Omega}[u] = I_{B_D}[u_0] +
\sum_{k=1}^{\infty} I_{B_D}[u_k]$.\ It suffices to prove that for
any\ $k=1,\ldots$,\ holds that
\[
I_{B_D}[u_k] \geq \frac{8(N-1)(N^2 -2N -2)}{(N^2 -4)^2} \int_{B_D}
|\Delta u_k|^2\, dx.
\]
Assume that the following inequality holds
\[
I_{B_D}[u_k] \geq a \int_{B_D} |\Delta u_k|^2\, dx,
\]
for some\ $0<a<1$\ and any\ $k =1, 2,...$.\ Taking into account
(\ref{3.Du}) and (\ref{3.I}) we obtain that
\begin{eqnarray*}
\int_{B_D} r^{2k-N+4} |\nabla g'|^2\, dx + \biggr ( 3 + 2k(N-3) +
\frac{N(N-4)}{2} \biggr ) \int_{B_D} r^{2k-N+2} |\nabla
g|^2\, dx \geq \;\;\;\;\;\;\;\;\;\;\;\;\;\;\;\;\;\;\;\;\; \\
\geq \frac{1}{1-a} \biggl \{ a \biggl [ \biggl ( \frac{N(N-4)}{4}
\biggr )^2 + \frac{N(N-4)}{2} (c_k +k^2) \biggr ] -
\frac{N(N-4)}{2}(c_k +k^2) \biggr \} \int_{B_D} r^{2k-N} g^2\, dx.
\end{eqnarray*}
Using now (\ref{2.17}) and (\ref{2.18}) we deduce that\ $a \leq
G(k)$,\ where
\[
G(k)= \frac{k^2 \biggl ( 3 +2k(N-3) + (k+1)^2 + \frac{N(N-4)}{2}
\biggr ) + \frac{N(N-4)}{2} \biggl ( 2k^2 +k (N-2) \biggr )}
{\biggl ( \frac{N(N-4)}{4} \biggr )^2 + k^2 \biggl ( 3 +2k(N-3) +
(k+1)^2 + \frac{N(N-4)}{2} \biggr ) + \frac{N(N-4)}{2} \biggl (
2k^2 +k (N-2) \biggr )}.
\]
However,\ $G(k)$\ is an increasing function for\ $k>1$.\ Hence,\
$a=G(1)=\frac{8(N-1)(N^2 -2N -2)}{(N^2 -4)^2}$\ and the proof is
completed.\ $\blacksquare$
\begin{lemma} \label{l4.2}
Let\ $u_0 \in C_{0}^{\infty} ([0,D])$.\ Then, the following
inequality holds
\begin{equation}\label{4.5}
I_{B_D}[u_0] \geq c\, \biggl ( \int_{B_D} |u_0|^{\frac{2N}{N-4}}
X^{\frac{2N-4}{N-4}} (\frac{|x|}{D})\, dx  \biggr
)^{\frac{N-4}{N}},
\end{equation}
for some positive constant\ $c$.\
\end{lemma}
\emph{Proof}\, Assume that\ $D=1$.\ From (\ref{3.6}) we have that
\begin{eqnarray} \label{4.6}
I_{B_D}[u_0] &\geq& c \int_{B_1} |x|^{4-N} |\Delta u_0|^2\, dx = c
\int_{0}^{1} r^{3} \biggl ( u^{''}_{0} + \frac{N-1}{r} u^{'}_{0}
\biggr )^2\, dr \nonumber \\
&=& c \biggl [ \int_{0}^{1} r^{3} ( u^{''}_{0} )^2\, dr +
(N-1)(N-2) \int_{0}^{1} r ( u^{'}_{0} )^2\, dr \biggr ] \nonumber \\
&=& c \biggl [ \int_{B_1 (\mathbb{R}^4)} ( \nabla u'_{0} )^2\, dx
+ (N-1)(N-2) \int_{B_1 (\mathbb{R}^4)} \frac{(u'_{0} )^2}{|x|^2}\,
dx \biggr ].
\end{eqnarray}
Applying now the Hardy inequality we have that
\begin{equation} \label{4.7}
\int_{B_1 (\mathbb{R}^4)} ( \nabla u'_{0} )^2\, dx \geq \biggl (
\frac{4-2}{2} \biggr )^2 \int_{B_1 (\mathbb{R}^4)} \frac{(u'_{0}
)^2}{|x|^2}\, dx \geq c \int_{0}^{1} r ( u^{'}_{0} )^2\, dr.
\end{equation}
So, from (\ref{4.6}) and (\ref{4.7}) we obtain that
\begin{equation} \label{4.8}
I[u_0] \geq c \int_{0}^{1} r ( u^{'}_{0} )^2\, dr.
\end{equation}
Next, we consider the following inequality
\begin{equation}\label{4.9}
\int_{0}^{1} r ( u^{'}_{0} )^2\, dr \geq c \biggl ( |u|^{q} r^{-1}
X^{1+q/2}(r)\, dr \biggr )^{2/q}
\end{equation}
which is implied from \cite[Theorem 3, p. 44]{maz85} with\ $X(t) =
(-logt)^{-1}$,\ $d\nu = r \chi_{[0,1]} dr$\ and\ $d\mu = r^{-1}
X^{\alpha} \chi_{[0,1]} dr$.\ Setting now\ $q= \frac{2N}{N-4}$,\
$\alpha = \frac{2N-4}{N-4}$\ and taking into account (\ref{4.8})
we conclude that
\begin{equation}\label{4.10}
I[u] \geq c \biggl ( \int_{B_1} |u_0|^{\frac{2N}{N-4}}
X^{\frac{2N-4}{N-4}}\, dx \biggr )^{\frac{N-4}{N}}.
\end{equation}
Following the same arguments we may prove that (\ref{4.10}) holds
for any\ $B_D$,\ $D>0$.\  $\blacksquare$
\vspace{0.2cm} \\
Using now Lemma \ref{l4.2} we prove inequality (\ref{4.5}) for
every\ $u \in C_{0}^{\infty} (\Omega)$. \vspace{0.3cm} \\
\emph{Proof of (\ref{4.11})}\, From inequality (\ref{4.1}) we have
that
\begin{equation}\label{4.12}
I_{\Omega}[u] \geq I_{B_D}[u_0] + c \int_{B_D} |\Delta
(u-u_0)|^2\, dx.
\end{equation}
The Sobolev imbedding and the fact that\ $X$\ is a bounded
function imply that
\begin{eqnarray}\label{4.13}
\int_{B_D} |\Delta (u-u_0)|^2\, dx &\geq& c \biggl ( \int_{B_D}
|u-u_0|^{\frac{2N}{N-4}}\, dx \biggr
)^{\frac{N-4}{N}} \nonumber \\
&\geq& c \biggl ( \int_{B_D} |u-u_0|^{\frac{2N}{N-4}}
X^{\frac{2N-4}{N-4}} (\frac{|x|}{D})\, dx \biggr
)^{\frac{N-4}{N}}.
\end{eqnarray}
Then from (\ref{4.5}), (\ref{4.12}) and (\ref{4.13}) we conclude
that
\[
I_{\Omega}[u] \geq c\, \biggl ( \int_{B_D} |u|^{\frac{2N}{N-4}}
X^{\frac{2N-4}{N-4}} (\frac{|x|}{D})\, dx  \biggr
)^{\frac{N-4}{N}} = c\, \biggl ( \int_{\Omega}
|u|^{\frac{2N}{N-4}} X^{\frac{2N-4}{N-4}} (\frac{|x|}{D})\, dx
\biggr )^{\frac{N-4}{N}}.\;\;\;\; \blacksquare
\]
\begin{theorem} \label{p4.4}
Let\ $\Omega$\ be a bounded domain,\ $D = sup_{x \in \Omega} |x|$\
and\ $u \in C_{0}^{\infty} (\Omega)$.\ Then,
\begin{equation} \label{4.14}
\mathbb{I}_{\Omega}[u] \geq \mathbb{I}_{B_D}[u_0] +
\frac{4(N-1)(N^2 -4N -4)}{(N^2 -4)^2} \int_{B_D} |\Delta
(u-u_0)|^2\, dx.
\end{equation}
\end{theorem}
\emph{Proof}\, Using the fact that\ $\mathbb{I}_{\Omega}[u] =
\mathbb{I}_{B_D}[u_0] + \sum_{k=1}^{\infty}
\mathbb{I}_{B_D}[u_k]$,\ it suffices to prove that for any\
$k=1,\ldots$,\ holds that
\[
\mathbb{I}_{B_D}[u_k] \geq \frac{4(N-1)(N^2 -4N -4)}{(N^2 -4)^2}
\int_{B_D} |\Delta u_k|^2\, dx.
\]
The result follows from (\ref{3.Du}) and (\ref{3.II}) using
(\ref{2.17}) and (\ref{2.18}).\ $\blacksquare$
\begin{lemma} \label{l4.5}
Let\ $u_0 \in C_{0}^{\infty} ([0,D])$.\ Then the following
inequality holds
\begin{equation}\label{4.19}
\mathbb{I}_{B_D}[u_0] \geq c\, \biggl ( \int_{B_D} |\nabla
u_0|^{\frac{2N}{N-2}} X^{1+\frac{N}{N-2}} (\frac{|x|}{D})\, dx
\biggr )^{\frac{N-2}{N}},
\end{equation}
for some positive constant\ $c$.\
\end{lemma}
\emph{Proof}\, Assume that\ $D=1$.\ Making some simple
calculations we may obtain that
\begin{eqnarray*}
\mathbb{I}[u_0] &=& \int_{0}^{1} r^{N-1} (u''_0 +
\frac{N-1}{r}u'_0)^2\, dr - \frac{N^2}{4} \int_{0}^{1} r^{N-3}
(u'_0)^2\, dr \\ &=& \int_{B_1} (u''_0)^2\, dx - \biggl (
\frac{N-2}{2} \biggr )^2 \int_{B_1} (u'_0)^2\, dx \\ &=&
\int_{B_1} |\nabla w|^2\, dx - \biggl ( \frac{N-2}{2} \biggr )^2
\int_{B_1} \frac{w^2}{|x|^2}\, dx,
\end{eqnarray*}
where\ $w = u'_0$.\ Using now the following inequality (see
\cite[Theorem A]{ft01})
\[
\int_{B_1} |\nabla w|^2\, dx - \biggl ( \frac{N-2}{2} \biggr )^2
\int_{B_1} \frac{w^2}{|x|^2}\, dx \geq c\, \biggl ( \int_{\Omega}
|w|^{\frac{2N}{N-2}} X^{1+\frac{N}{N-2}}\, dx \biggr
)^{\frac{N-2}{N}}
\]
which hold for any\ $w \in H_{0}^{1}(B_1)$,\ we obtain that
(\ref{4.19}) holds for any\ $u(r) \in C_{0}^{\infty} (B_1)$.\
Then, following the same arguments we may obtain that (\ref{4.19})
hold for any\ $B_D$,\ $D>0$.\ $\blacksquare$ \vspace{0.3cm} \\
\emph{Proof of (\ref{4.20})}\, As in the proof of (\ref{4.11}) the
result is a consequence of Proposition \ref{p4.4}, Lemma
\ref{l4.5} and of the following inequality
\begin{eqnarray*}
\int_{B_D} |\Delta (u-u_0)|^2\, dx &\geq& c \biggl ( \int_{B_D}
|\nabla u- \nabla u_0|^{\frac{2N}{N-2}}\, dx \biggr
)^{\frac{N-2}{N}} \nonumber \\
&\geq& c \biggl ( \int_{B_D} |\nabla u- \nabla
u_0|^{\frac{2N}{N-2}} X^{1+\frac{N}{N-2}} (\frac{|x|}{D})\, dx
\biggr )^{\frac{N-2}{N}},
\end{eqnarray*}
which is implied from the Sobolev imbedding and the fact that\
$X$\ is a bounded function.\ $\blacksquare$ \vspace{0.2cm} \\
\emph{Proof of Theorem \ref{t5.4}}\, Inequality (\ref{5.8}) is an
immediate consequence from Proposition  \ref{p3.1} and Inequality
(\ref{1.3}).\ The fact that\ $\biggl ( 1 + \frac{N(N-4)}{8} \biggr
)$\ is the best constant will be establish in Section 4.
$\blacksquare$ \vspace{0.3cm} \\
\emph{Proof of Theorem \ref{t5.5}}\, We decompose\ $u$\ into
spherical harmonics. Then, the result is an immediate consequence
of (\ref{1.3}), (\ref{3.Du}) and (\ref{2.12}), with\ $V(|x|) =
\frac{1}{4} (N^2 + X_{1}^{2}(\frac{|x|}{D}) \ldots
X_{i}^{2}(\frac{|x|}{D}) )$.\ The fact that\ $\frac{N^2}{4}$\ and\
$\frac{1}{4}$\ are the best constants will be establish in Section
4. $\blacksquare$
%
%
%
%
\section{Best Constants}
\setcounter{equation}{0}
Throughout this section we may assume that\ $\Omega$\ is a bounded
domain, such that\ $B_1(0) \subset \Omega$ and\ $N \geq 5$.\ We
initially establish that the constants appearing in the
inequalities of Section 2.2 are the best ones. For some\ $\epsilon
>0$\ and\ $0<a_1$\ we introduce the minimizing sequences\
$u^{\epsilon}$\ and\ $v^{\epsilon}$\ to be defined as:
\[
u^{\epsilon} := r^{-\frac{N-4}{2}+\epsilon}
X_{1}^{\frac{-1+a_1}{2}} \phi(r),\;\;\;\;\;\; v^{\epsilon} :=
r^{\frac{N-4}{2}} u^{\epsilon} = r^{\epsilon}
X_{1}^{\frac{-1+a_1}{2}} \phi(r),
\]
where\ $X_1(t) = (1- \log t)^{-1}$\ and\ $\phi(r) \in
C^{\infty}_{0} (B_1)$\ is a smooth cutoff function, such that\ $0
\leq \phi \leq 1$, with\ $\phi \equiv 1$\ in\ $B_{1/2}$.\
\begin{lemma} \label{l5.1}
As\ $\epsilon \to 0^+$\ and\ $a_1 \to 0^+$,\ we have
\begin{eqnarray}
\label{5.1a} i)\;\;  \frac{1}{c_N}\;  \int_{\Omega} |x|^{2-N}
|\nabla v^{\epsilon}|^2\, dx = - \frac{-1+a_{1}}{4} \int_{0}^{1}
r^{-1 + 2 \epsilon}\, X_{1}^{1+a_1}\, \phi^2(r)\, dr + O(1), \\
\label{5.1b} ii)\;\;  \frac{1}{c_N}\;  \int_{\Omega} |x|^{4-N}
|\Delta v^{\epsilon}|^2\, dx = - \frac{-1+a_1}{4} (N-2)^2
\int_{0}^{1} r^{2\epsilon-1} X_{1}^{1+a_1} \phi^2\, dr + O(1),\\
\label{5.1c} iii)\;\;  \frac{1}{c_N}\;  \int_{\Omega} |x|^{-2}
|\nabla u^{\epsilon}|^2\, dx = - \frac{-1+a_{1}}{4} \int_{0}^{1}
r^{2\epsilon-1} X_{1}^{1+a_1} \phi^2\, dr \nonumber
\\  + \biggl ( \frac{N-4}{2} \biggr )^2 \int_{0}^{1}
r^{2\epsilon-1} X_{1}^{-1+a_{1}} \phi^2\, dr + O(1), \\
\label{5.1d} iv)\;\; \frac{1}{c_N}\;  \int_{\Omega} |\Delta
u^{\epsilon}|^2\, dx = - \frac{-1+a_{1}}{8} (N^2-4N+8)
\int_{0}^{1} r^{2\epsilon-1} X_{1}^{1+a_1} \phi^2\, dr \nonumber
\\ + \biggl ( \frac{N(N-4)}{4} \biggr )^2 \int_{0}^{1}
r^{2\epsilon-1}
X_{1}^{-1+a_{1}} \phi^2\, dr +O(1), \\
%
\\
\label{5.1e} v)\;\;  \frac{1}{c_N}\; I[u^{\epsilon}] = -
\frac{-1+a_{1}}{8} (N^2-4N+8) \int_{0}^{1} r^{2\epsilon-1}
X_{1}^{1+a_1} \phi^2\, dr +O(1),
\\
\label{5.1f} vi)\;\;  \frac{1}{c_N}\; \mathbb{I}[u^{\epsilon}] = -
\frac{-1+a_{1}}{16} (N-4)^2 \int_{0}^{1} r^{2\epsilon-1}
X_{1}^{1+a_1} \phi^2\, dr +O(1),
\end{eqnarray}
where\ $c_N$\ is the volume of the unit sphere in\
$\mathbb{R}^N$.\
\end{lemma}
\emph{Proof}\, The conclusion follows from the properties of the
functions\ $X_1$,\ $\phi$\ and standard arguments based on
integration by parts which also imply that
\[
\epsilon \int_{0}^{1} r^{-1 + 2 \epsilon}\, X_{1}^{-1+a_1}\,
\phi^2(r)\, dr = - \frac{-1+a_{1}}{2} \int_{0}^{1} r^{-1 + 2
\epsilon}\, X_{1}^{a_1}\, \phi^2(r)\, dr + O(1)
\]
and
\[
2 \epsilon \int_{0}^{1} r^{-1 + 2 \epsilon}\, X_{1}^{a_1}\,
\phi^2(r)\, dr = - a_1 \int_{0}^{1} r^{-1 + 2 \epsilon}\,
X_{1}^{1+a_1}\, \phi^2(r)\, dr + O(1).
\]
\setcounter{theorem}{1}
\begin{theorem} \label{p5.1}
The quantities
\begin{eqnarray*}
&\mbox{i)}&  4 + \frac{N(N-4)}{2}\;\;\;\;\;\;  \mbox{in
inequality (\ref{3.1})}, \\
&\mbox{ii)}&  \frac{1}{2} + \frac{2}{(N-2)^2}\;\;\;\;\;\;
\mbox{in inequality (\ref{3.6})}, \\
&\mbox{iii)}& \biggl ( \frac{N-4}{2} \biggr )^2\;\;\;\;\;\;
\mbox{in inequality (\ref{3.7})}, \\
&\mbox{iv)}& 2(N-2)^2\;\;\;\;\;\; \mbox{in inequality
(\ref{3.3})},
\\
&\mbox{v)}& \biggl ( \frac{N-4}{2(N-2)} \biggr )^2\;\;\;\;\;\;
\mbox{in inequality (\ref{3.8})}, \\
%
%
&\mbox{vi)}& \frac{N^2}{4}\;\;\;\;\;\; \mbox{in inequality
(\ref{rel2})},
\end{eqnarray*}
are the best constants.
\end{theorem}
\emph{Proof}\, i)\ Relations (\ref{5.1a}) and (\ref{5.1e}) imply
that
\begin{eqnarray*}
\frac{I[u^{\epsilon}]}{\int_{\Omega} |x|^{2-N} |\nabla
v^{\epsilon}|^2\, dx} &=& \frac{\frac{-1+a_1}{8} (N^2-4N+8)
\int_{0}^{1} r^{2\epsilon-1} X_{1}^{1+a_1} \phi^2\, dr +O(1)}{-
\frac{-1+a_1}{4} \int_{0}^{1} r^{-1 + 2 \epsilon}\,
X_{1}^{1+a_1}\, \phi^2(r)\, dr +
O(1)} \\
&\to& 4 + \frac{N(N-4)}{2},
\end{eqnarray*}
as\ $\epsilon \downarrow 0$\ and\ $a_1 \downarrow 0$.\ In the same
way the conclusion follows for the cases ii) - v). For the last
case, observe that
\[
\frac{\int_{0}^{1} r^{2\epsilon-1} X_{1}^{-1+a_1} \phi^2\,
dr}{\int_{0}^{1} r^{2\epsilon-1} X_{1}^{1+a_{1}} \phi^2\, dr} \to
\infty,
\]
as\ $\epsilon \downarrow 0$\ and\ $ a \downarrow 0$.\ Then, from
(\ref{5.1c}) and (\ref{5.1d}) we derive that
\begin{eqnarray*}
\frac{\int_{\Omega} |\Delta u^{\epsilon}|^2\, dx}{\int_{\Omega}
\frac{|\nabla u^{\epsilon}|^2}{|x|^2}\, dx} &=&
\frac{\frac{-1+a_{1}}{8} (N^2-4N+8) \int_{0}^{1} r^{2\epsilon-1}
X_{1}^{1+a_1} \phi^2\, dr + \biggl ( \frac{N(N-4)}{4} \biggr )^2
\int_{0}^{1} r^{2\epsilon-1} X_{1}^{-1+a_{1}} \phi^2\, dr
}{\frac{-1+a_{1}}{4} \int_{0}^{1} r^{2\epsilon-1} X_{1}^{1+a_{1}}
\phi^2\, dr + \biggl ( \frac{N-4}{2} \biggr )^2
\int_{0}^{1} r^{2\epsilon-1} X_{1}^{-1+a_{1}} \phi^2\, dr } \\
&\to& \frac{N^2}{4},
\end{eqnarray*}
as\ $\epsilon \downarrow 0$\ and\ $a_{1} \downarrow 0$.\
$\blacksquare$\ \vspace{0.2cm}

Next we complete the proofs of Theorems \ref{t5.4} and \ref{t5.5}.
We introduce the minimizing sequences for the k-Improved
Hardy-Rellich inequalities. For small positive parameters\
$\epsilon,\, a_1,\, a_2,\,...,a_k $\ we define
\[
u(x) := w(x)\, \phi(|x|),\;\;\;
w(x) := |x|^{-\frac{N-4}{2}+\epsilon} X_{1}^{\frac{-1+a_1}{2}}
X_{2}^{\frac{-1+a_2}{2}} \cdots X_{k}^{\frac{-1+a_k}{2}},
\]
where\ $\phi$\ is the previous test function and\ $X_m = X_1
(X_{m-1})$,\ $m=2,...,k$.\ To prove the results we shall estimate
the corresponding Rayleigh quotients of $u$ (\ref{ray1}),
(\ref{ray2}) in the limit\ $\epsilon \to 0$,\ $a_1 \to 0$,\
$...$,\ $a_k \to 0$\ in this order.

In the sequel we shall repeatedly use the differentiant rule
\[
\frac{d}{dt} X_{i}^{\beta}(t) = \frac{\beta}{t} X_1 X_2 \cdots
X_{i-1} X_{i}^{1+\beta},\;\;\; \beta \ne -1,\; i=1,2,...,
\]
and with integrals of the form
\[
Q = \int_{0}^{1} r^{-1+2\epsilon} X_{1}^{1+\beta_1}
X_{2}^{1+\beta_2} \cdots X_{k}^{1+\beta_k} \phi^2 (r) \, dr.
\]
For this we notice that
\[
Q < \infty \Leftrightarrow
\left \{
\begin{array}{ll}
\epsilon>0,\;\;\;\;    or, \\
\epsilon = 0\;
\mbox{and}\; \beta_1 >0,\;\;\;\;   or, \\
\epsilon = 0,\, \beta_1 =0\;
\mbox{and}\; \beta_1 >0,\;\;\;   or, \\
\vdots \\
\epsilon = 0,\, \beta_1 =0,\, ...,\, \beta_{k-1}=0\;
\mbox{and}\;\; \beta_k
>0.
\end{array}
\right.
\]
Also as we pass to the limit\ \ $\epsilon \to 0$,\ $a_1 \to 0$,\
$...$,\ $a_k \to 0$\ we have
\begin{eqnarray*}
\int_{\Omega} |\Delta u|^2\, dx &=& \int_{\Omega} |\Delta w|^2
\phi^2\, dx + O(1), \\
\int_{\Omega} \frac{|\nabla u|^2}{|x|^2}\, dx &=& \int_{\Omega}
\frac{|\nabla w|^2}{|x|^2} \phi^2\, dx + O(1), \\
\int_{\Omega} \frac{|\nabla u|^2}{|x|^2}\, X_{1}^{2} \cdots
X_{i}^{2} dx &=& \int_{\Omega} \frac{|\nabla w|^2}{|x|^2}
X_{1}^{2} \cdots X_{i}^{2} \phi^2\, dx + O(1),\;\;\; i=1,...,k.
\end{eqnarray*}
It is not difficult to see that
\[
\nabla w(x) = |x|^{-\frac{N-2}{2}+\epsilon}
X_{1}^{\frac{-1+a_1}{2}} X_{2}^{\frac{-1+a_2}{2}} \cdots
X_{k}^{\frac{-1+a_k}{2}} \biggl [- \frac{N-4}{2} + \epsilon +
\frac{1}{2} \eta(x) \biggr ] \frac{x}{|x|},
\]
where\ $\eta(x) = (-1+a_1)X_1 + (-1+a_2)X_1 X_2 + ... +
(-1+a_k)X_1 \cdots X_k$\ and
\[
\Delta w(x) = \frac{1}{4} |x|^{-\frac{N}{2}+\epsilon}
X_{1}^{\frac{-1+a_1}{2}} X_{2}^{\frac{-1+a_2}{2}} \cdots
X_{k}^{\frac{-1+a_k}{2}} \biggl [ -N(N-4)+ 8 \epsilon +4
\epsilon^2 +4(1+\epsilon) \eta(x) + \eta^2 (x) + 2 B(x) \biggr ],
\]
where
\[
\begin{array}{ccll}
B(|x|) &=& (-1+a_1)X_{1}^{2} + (-1+a_2) (X_{1}^{2} X_2 + X_{1}^{2}
X_{2}^{2}) + \ldots \\
&& \;\;\;\;\;\;\;\;\; +(-1+a_k) (X_{1}^{2} X_{2} \cdots X_{k} +
\ldots + X_{1}^{2} X_{2}^{2} \cdots X_{k}^{2}) \\
&=& \sum_{i=1}^{k} (-1+a_i) X_{1}^{2}\cdots X_{i}^{2} +
\sum_{i=2}^{k} \sum_{j=1}^{i-1} (-1+a_i) X_{1}^{2}\cdots X_{j}^{2}
X_{j+1} \cdots X_{i} \\
&=& \sum_{i=1}^{k} (-1+a_i) X_{1}^{2}\cdots X_{i}^{2} +
\sum_{j=1}^{k-1} \sum_{i=j+1}^{k} (-1+a_i) X_{1}^{2}\cdots
X_{j}^{2} X_{j+1} \cdots X_{i},
\end{array}
\]
Note also that
\[
r \eta'(r) = B(r)
\]
and
\[
\eta^2 (x) = \sum_{i=1}^{k} (-1+a_i)^2 X_{1}^{2} \cdots X_{i}^{2}
+ 2 \sum_{j=1}^{k-1} \sum_{i=j+1}^{k} (-1+a_i)(-1+a_j) X_{1}^{2}
\cdots X_{j}^{2} X_{j+1} \cdots X_{i}.
\]
Then we have that
\[
\begin{array}{ccll}
\int_{\Omega} |\Delta u|^2 dx &=& \int_{\Omega} \frac{w^2}{|x|^4}
\biggl [ \biggl ( - \frac{N(N-4)}{4} +2 \epsilon + \epsilon^2
\biggr )^2 + (1+\epsilon)^2 \eta^2 + (\frac{1}{4} \eta^2 +
\frac{1}{2}B)^2\\
&& + 2 (1+\epsilon) \biggl ( - \frac{N(N-4)}{4} +2 \epsilon +
\epsilon^2 \biggr ) \eta + 2 \biggl (-\frac{N(N-4)}{4}+2 \epsilon
+ \epsilon^2 \biggr ) (\frac{1}{4} \eta^2 + \frac{1}{2}B)\\
&&  +2 (1+\epsilon)(\frac{1}{4}
\eta^2 + \frac{1}{2}B)\eta \biggr ] \phi^2\, dx \\
&=& \int_{\Omega} \frac{w^2}{|x|^4} \biggl [ \biggl ( -
\frac{N(N-4)}{4} +2 \epsilon + \epsilon^2 \biggr )^2 +
(1+\epsilon)^2 \eta^2 \\
&& + 2 (1+\epsilon) \biggl ( - \frac{N(N-4)}{4} +2 \epsilon +
\epsilon^2 \biggr ) \eta + 2 \biggl (-\frac{N(N-4)}{4}+2 \epsilon
+ \epsilon^2 \biggr ) (\frac{1}{4} \eta^2 + \frac{1}{2}B) \biggr ]
\cdot \\ && \cdot \phi^2\, dx +O(1),
\end{array}
\]
\[
\begin{array}{ccll}
\int_{\Omega} \frac{|\nabla u|^2}{|x|^2} X_{1}^{2} \cdots
X_{i}^{2}\, dx = \int_{\Omega} \frac{w^2}{|x|^4} \biggl [ \biggl (
-\frac{N-4}{2} + \epsilon \biggr )^2 + \biggl ( -\frac{N-4}{2} +
\epsilon \biggr ) \eta + \frac{1}{4} \eta^2 \biggr ]\, X_{1}^{2}
\cdots X_{i}^{2}\, \phi^2\, dx +O(1), \vspace{0.2cm} \\
\int_{\Omega} \frac{u^2}{|x|^4} X_{1}^{2} \cdots X_{i}^{2}\, dx =
\int_{\Omega} \frac{w^2}{|x|^4} X_{1}^{2} \cdots X_{i}^{2}\,
\phi^2 dx +O(1).
\end{array}
\]
An important quantity that appears is
\[
\sum_{i=1}^{k} a_i A_i  - \sum_{i=1}^{k-1} \sum_{j=i+1}^{k}
(1-a_j) \Gamma_{ij},
\]
where
\begin{eqnarray*}
A_i (a_1, \ldots a_k) &:=& \int_{0}^{1} r^{-1} X_{1}^{1+a_1}
\cdots X_{i}^{1+a_i}
X_{i+1}^{-1+a_{i+1}} \cdots X_{k}^{-1+a_k} \phi^2\, dr, \;\;\; i=1, \ldots, k \\
\Gamma_{ij} (a_1, \ldots a_k) &:=& \int_{0}^{1} r^{-1}
X_{1}^{1+a_1} \cdots X_{i}^{1+a_i} X_{i+1}^{a_{i+1}} \cdots
X_{j}^{a_j} X_{j+1}^{-1+a_{j+1}} \cdots X_{k}^{-1+a_k} \phi^2\,
dr,\;\;\; i<j.
\end{eqnarray*}
We will pass to the limit initially\ $a_1 \to 0^+$\ and then\ $a_2
\to 0^+,\ ...,\ a_{k-1} \to 0^+$. In passing to the limit we will
use identities similar to the ones used in Step 8 of \cite{bft03},
in particular we have
\[
a_1 A_1 = \int_{0}^{1} (X_{1}^{a_1})' X_{2}^{-1+a_2} \cdots
X_{k}^{-1+a_k} \phi^2\, dr =  - \sum_{j=2}^{k} (-1+a_j)
\Gamma_{ij} + O(1),\;\;\; \mbox{as}\;\; a_1 \downarrow 0.
\]
Therefore,
\begin{equation}\label{eqb3}
\sum_{i=1}^{k} a_i A_i - \sum_{i=1}^{k-1} \sum_{j=i+1}^{k} (1-a_j)
\Gamma_{ij} = \sum_{i=2}^{k} a_i A_i - \sum_{i=2}^{k-1}
\sum_{j=i+1}^{k} (1-a_j) \Gamma_{ij} + O(1),\;\;\; \mbox{as}\;\;
a_1 \downarrow 0
\end{equation}
Then we pass to the limit\ $a_1 \to 0$,\ in the right hand side of
(\ref{eqb3}). Again we use the identity
\begin{eqnarray*}
a_2 A_2 (0,a_2,...,a_k) &=& \int_{0}^{1} (X_{2}^{a_2})'
X_{3}^{-1+a_3} \cdots X_{k}^{-1+a_k} \phi^2\, dr \\ &=& -
\sum_{j=3}^{k} (-1+a_j) \Gamma_{ij} + O(1),\;\;\; \mbox{as}\;\;
a_2 \downarrow 0.
\end{eqnarray*}
and therefore by iterating the previous procedure we pass to the
limit $a_2 \to 0^+,\ ...,\ a_{k-1} \to 0^+$ to conclude
\begin{eqnarray} \label{AG}
\sum_{i=1}^{k} a_i A_i  - \sum_{i=1}^{k-1} \sum_{j=i+1}^{k}
(1-a_j) \Gamma_{ij} &=& a_k A_k(0,0,\ldots, 0, a_k) +O(1), \nonumber \\
&=& a_k \int_{0}^{1} r^{-1} X_1 X_2 \cdots X_{k-1} X_{k}^{-1+a_k}
\phi^2\, dr +O(1), \mbox{as}\; a_k \downarrow 0.
\end{eqnarray}
\emph{Completion of Proof of Theorem \ref{t5.4}}. We use the
previous test functions to conclude that
\[
\begin{array}{ccclll}
R[u] := \int_{\Omega} |\Delta u|^2 dx - \biggl ( \frac{N(N-4)}{4}
\biggr )^2 \int_{\Omega} \frac{u^2}{|x|^4} dx - \biggl (1+
\frac{N(N-4)}{8} \biggr ) \sum_{i=1}^{k-1} \int_{\Omega}
\frac{u^2}{|x|^4} X_{1}^{2} \cdots
X_{i}^{2}\, dx  = \\
= \int_{\Omega} \frac{w^2}{|x|^4} \biggl [ \epsilon^2
(2+\epsilon)^2 - \frac{N(N-4)}{2} \epsilon (2+\epsilon) +2
(1+\epsilon) \biggl (
- \frac{N(N-4)}{4} +2 \epsilon + \epsilon^2 \biggr ) \eta \\
+ \biggl (1 - \frac{N(N-4)}{8} + 3 \epsilon + \frac{3}{2}
\epsilon^2 \biggr ) \eta^2 + \biggl ( - \frac{N(N-4)}{4} +2
\epsilon + \epsilon^2 \biggr ) B \\
- \biggl (1+ \frac{N(N-4)}{8} \biggr ) \sum_{i=1}^{k-1} X_{1}^{2}
\cdots X_{i}^{2} \biggr ]\, \phi^2\, dx +O(1), \\
= c_N \int_{0}^{1} r^{-1+2\epsilon} X_{1}^{-1+a_1} \cdots
X_{k}^{-1+a_k} \biggl [ \epsilon^2 (2+\epsilon)^2 -
\frac{N(N-4)}{2} \epsilon
(2+\epsilon)\;\;\;\;\;\;\;\;\;\;\;\;\;\;\;
\;\;\;\;\;\; \\
+2 (1+\epsilon) \biggl ( - \frac{N(N-4)}{4} +2 \epsilon +
\epsilon^2 \biggr ) \eta + \biggl (1 - \frac{N(N-4)}{8} + 3
\epsilon + \frac{3}{2} \epsilon^2 \biggr ) \eta^2 \\
+ \biggl ( - \frac{N(N-4)}{4} +2 \epsilon + \epsilon^2 \biggr ) B
- \biggl (1+ \frac{N(N-4)}{8} \biggr ) \sum_{i=1}^{k-1} X_{1}^{2}
\cdots X_{i}^{2} \biggr ]\, \phi^2\, dr +O(1),
\end{array}
\]
But
\begin{eqnarray} \label{eqb1}
2 \epsilon \int_{0}^{1} r^{-1+2\epsilon} X_{1}^{-1+a_1} \cdots
X_{k}^{-1+a_k} \phi^2(r) dr &=& - \int_{0}^{1} r^{2\epsilon}
(X_{1}^{-1+a_1} \cdots X_{k}^{-1+a_k})' \phi^2(r) dr +O(1) \nonumber \\
&=& - \int_{0}^{1} r^{-1+2\epsilon} X_{1}^{-1+a_1} \cdots
X_{k}^{-1+a_k} \eta \phi^2(r) dr +O(1).\;\;\;\;
\end{eqnarray}
and
\begin{eqnarray} \label{eqb2}
2 \epsilon \int_{0}^{1} r^{-1+2\epsilon} X_{1}^{-1+a_1} \cdots
X_{k}^{-1+a_k} \eta \phi^2(r) dr = - \int_{0}^{1} r^{-1+2\epsilon}
X_{1}^{-1+a_1} \cdots X_{k}^{-1+a_k} \eta^2 \phi^2(r) dr \nonumber
\\ - \int_{0}^{1} r^{-1+2\epsilon} X_{1}^{-1+a_1} \cdots
X_{k}^{-1+a_k} B \phi^2(r) dr +O(1).
\end{eqnarray}
Therefore
\[
\begin{array}{ccll}
R[u] = c_N \int_{0}^{1} r^{-1+2\epsilon} X_{1}^{-1+a_1} \cdots
X_{k}^{-1+a_k} \biggl [ \epsilon^3 +\epsilon^4 +(6+2\epsilon)
\epsilon^2 \eta
+ \biggl (3 \epsilon + \frac{3}{2} \epsilon^2 \biggr ) \eta^2
\;\;\;\;\;\;\;\;\;\;\;\;\;\;\;\;\;\;\;\;\;\;\;\; \\
+ \biggl (-1 - \frac{N(N-4)}{8} +2 \epsilon + \epsilon^2 \biggr )
B
- \biggl (1+ \frac{N(N-4)}{8} \biggr ) \sum_{i=1}^{k-1} X_{1}^{2}
\cdots X_{i}^{2} \biggr ]\, \phi^2\, dr +O(1),
\end{array}
\]
passing to the limit $\epsilon \to 0$, we obtain
\[
\begin{array}{ccll}
 \frac{1}{c_N} R[u] &=& \biggl(1+\frac{N(N-4)}{8}
\biggr) \int_{0}^{1} r^{-1} X_{1}^{-1+a_1} \cdots X_{k}^{-1+a_k}
\biggl [ B + \sum_{i=1}^{k-1}
X_{1}^{2} \cdots X_{i}^{2} \biggr ]\, \phi^2\, dr +O(1) \\
&=& \biggl(1+\frac{N(N-4)}{8} \biggr) \int_{0}^{1} r^{-1}
X_{1}^{1+a_1} \cdots X_{k}^{1+a_k} \phi^2\, dr - \\
&& - \biggl(1+\frac{N(N-4)}{8} \biggr) \int_{0}^{1} r^{-1}
X_{1}^{-1+a_1} \cdots X_{k}^{-1+a_k} \biggl [ \sum_{i=1}^{k} a_i
X_{1}^{2} \cdots X_{i}^{2} + \\
&& + \sum_{j=1}^{k-1} \sum_{i=j+11}^{k} (-1+a_i) X_{1}^{2} \cdots
X_{j}^{2} X_{j+1} \cdots X_{i} \biggr ]\, \phi^2\, dr +O(1).
\end{array}
\]
or
\begin{eqnarray*}
\frac{1}{c_N} R[u] = \biggl(1+\frac{N(N-4)}{8} \biggr) A_k -
\biggl(1+\frac{N(N-4)}{8} \biggr) \biggl( \sum_{i=1}^{k} a_i A_i -
\sum_{i=1}^{k-1}\sum_{j=i+1}^{k} (1-a_j)\Gamma_{ij} \biggr) +O(1).
\end{eqnarray*}
However, we can pass to the limit\ $a_1 \downarrow 0,...a_{k-1}
\downarrow 0$\ see (\ref{AG}), to conclude that
\[
\frac{1}{c_N} R[u] = \biggl(1+\frac{N(N-4)}{8} \biggr) A_k -
\biggl(1+\frac{N(N-4)}{8} \biggr) a_k A_k + O(1),\; \mbox{as}\;
a_k \downarrow 0.
\]
The Rayleigh quotient now of (\ref{ray1}) is smaller or equal than
\[
\frac{\biggl(1+\frac{N(N-4)}{8} \biggr) A_k -
\biggl(1+\frac{N(N-4)}{8} \biggr) a_k A_k +O(1)}{A_k} \to
1+\frac{N(N-4)}{8},
\]
since\ $A_k \to \infty$,\ as\ $a_k \downarrow 0$.\ $\blacksquare$
\vspace{0.2cm} \\
\emph{Completion of Proof of Theorem \ref{t5.5}}. Once more we use
the same minimizing sequence to conclude
\[
\begin{array}{ccll}
\int_{\Omega} |\Delta u|^2 dx - \frac{N^2}{4} \int_{\Omega}
\frac{|\nabla u|^2}{|x|^2} dx - \frac{1}{4} \sum_{i=1}^{k-1}
\int_{\Omega} \frac{|\nabla u|^2}{|x|^2} X_{1}^{2} \cdots
X_{i}^{2}\, dx
=\;\;\;\;\;\;\;\;\;\;\;\;\;\;\;\;\;\;\;\;\;\;\;\;\;\;\;\;\;
\;\;\;\;\;\;\;\;\;\;\;\;\; \\
= c_N \int_{0}^{1} r^{-1+2\epsilon} X_{1}^{-1+a_1} \cdots
X_{i}^{-1+a_k} \biggl [ \biggl ( -\frac{N(N-4)}{4} +2 \epsilon +
\epsilon^2 \biggr )^2
+(1+\epsilon)^2 \eta^2 \\
+2(1+\epsilon) \biggl ( -\frac{N(N-4)}{4} +2 \epsilon + \epsilon^2
\biggr ) \eta +2 \biggl ( -\frac{N(N-4)}{4} +2 \epsilon +
\epsilon^2 \biggr ) (\frac{1}{4} \eta^2 +\frac{1}{2} B)  \\
-\frac{N^2}{4} \biggl ( -\frac{N-4}{2} + \epsilon \biggr )^2
-\frac{N^2}{4} \biggl ( -\frac{N-4}{2} + \epsilon \biggr ) \eta -
\frac{N^2}{16} \eta^2 \\
-\frac{1}{4} \biggl ( -\frac{N-4}{2} + \epsilon \biggr )^2
\sum_{i=1}^{k-1} X_{1}^{2} \cdots X_{i}^{2} \biggr ]\, \phi^2\, dr
+O(1), \\
= c_N \int_{0}^{1} r^{-1+2\epsilon} X_{1}^{-1+2a_1} \cdots
X_{i}^{-1+2a_k} \biggl [ \frac{N(N-4)^2}{4}\epsilon + \biggl ( 4 -
\frac{N(N-4)}{2} - \frac{N^2}{4} \biggr ) \epsilon^2 +4 \epsilon^3
+ \epsilon^4 \\
+ \biggl ( \frac{N(N-4)^2}{8} + \biggl ( 4 - \frac{N(N-4)}{2} -
\frac{N^2}{4} \biggr ) \epsilon +6 \epsilon^2 +2\epsilon^3 \biggr
) \eta \\
+ \biggl ( 1- \frac{N(N-4)}{8} -\frac{N^2}{16} +3\epsilon +
\frac{3}{2} \epsilon^2 \biggr ) \eta^2 + \biggl ( -
\frac{N(N-4)}{4} +2\epsilon + \epsilon^2 \biggr ) B \\
- \frac{1}{4} \biggl ( -\frac{N-4}{2} + \epsilon \biggr )^2
\sum_{i=1}^{k-1} X_{1}^{2} \cdots X_{i}^{2} \biggr ]\, \phi^2\, dr
+O(1),
\end{array}
\]
We now use identities (\ref{eqb1}), (\ref{eqb2}) and passing to
the limit\ $\epsilon \to 0$, to conclude that
\begin{eqnarray*}
\int_{\Omega} |\Delta u|^2 dx - \frac{N^2}{4} \int_{\Omega}
\frac{|\nabla u|^2}{|x|^2} dx - \frac{1}{4} \sum_{i=1}^{k-1}
\int_{\Omega} \frac{|\nabla u|^2}{|x|^2} X_{1}^{2} \cdots
X_{i}^{2}\, dx
=\;\;\;\;\;\;\;\;\;\;\;\;\;\;\;\;\;\;\;\;\;\;\;\;\;\;\;\;\;
\;\;\;\;\;\;\;\;\;\;\;\;\; \\
= - \frac{1}{4} \biggl ( \frac{N-4}{2} \biggr )^2 c_N \int_{0}^{1}
r^{-1} X_{1}^{-1+a_1} \cdots X_{i}^{-1+a_k} \biggl [ B +
\sum_{i=1}^{k-1} X_{1}^{2} \cdots X_{i}^{2} \biggr ]\, \phi^2\, dr
+O(1), \\
= \frac{1}{4} \biggl ( \frac{N-4}{2} \biggr )^2 c_N A_k -
\frac{1}{4} \biggl ( \frac{N-4}{2} \biggr )^2 c_N \biggl(
\sum_{i=1}^{k} a_i A_i - \sum_{i=1}^{k-1}\sum_{j=i+1}^{k}
(1-a_j)\Gamma_{ij} \biggr) +O(1).
\end{eqnarray*}
As before, using (\ref{AG}) we can pass to the limit\ $a_1
\downarrow 0,...a_{k-1} \downarrow 0$\ see (\ref{AG}), the
Rayleigh quotient now of (\ref{ray2}) is smaller or equal than
\[
\frac{ \frac{1}{4} \biggl ( \frac{N-4}{2} \biggr )^2 A_k -
\frac{1}{4} \biggl ( \frac{N-4}{2} \biggr )^2 a_k A_k +O(1)}{
\biggl ( \frac{N-4}{2} \biggr )^2 A_k +O(1)} \to \frac{1}{4},
\]
since\ $A_k \to \infty$,\ as\ $a_k \downarrow 0$.\ $\blacksquare$
%
%
\section{Existence of minimizers in\
$W^{2,2}_0 (\Omega, |x|^{-(N-4)})$}
\setcounter{equation}{0}
In this section we assume certain improved inequalities
(\ref{rel1}), (\ref{rel2}) in\ $v$-terms and we prove the
existence of minimizers, in some appropriate weighted spaces.
Assume that\ $\Omega \subset \mathbb{R}^N$\ is a bounded domain
containing the origin and\ $N \geq 5$.\ We introduce the space\
$W^{2,2}_0 (\Omega, |x|^{-(N-4)})$\ to be defined as the closure
of the\ $C_{0}^{\infty}$\ functions with respect to the norm
\begin{eqnarray} \label{6.1}
||v||^{2}_{W} &:=& \int_{\Omega} |x|^{-(N-4)} |\Delta v|^2\, dx +
\int_{\Omega} |x|^{-N} |x \cdot \nabla v|^2\, dx + \int_{\Omega}
|x|^{-(N-2)} |\nabla v|^2\, dx \nonumber \\ && + \int_{\Omega}
|x|^{-N}\, v^2\, dx.
\end{eqnarray}
It is clear that\ $W^{2,2}_0 (\Omega, |x|^{-(N-4)})$\ is a Hilbert
space with inner product
\begin{eqnarray*}
< \phi, \psi >_{W} &:=& \int_{\Omega} |x|^{-(N-4)} \Delta \phi\,
\Delta \psi\, dx + \int_{\Omega} |x|^{-N} (x \cdot \nabla \phi) (x
\cdot \nabla \psi)\, dx \\ && + \int_{\Omega} |x|^{-(N-2)} \nabla
\phi \cdot \nabla \psi\, dx + \int_{\Omega} |x|^{-N} \phi\, \psi\,
dx.
\end{eqnarray*}
\begin{lemma} \label{l6.1}

i)\ If\ $u \in H_{0}^{2} (\Omega)$,\ then\ $v =
|x|^{\frac{N-4}{2}} u \in W^{2,2}_0 (\Omega, |x|^{-(N-4)})$.\

ii) If\ $v \in C_{0}^{\infty} (\Omega)$,\ then\ $u = |x|^{-\alpha}
v \in H_{0}^{2} (\Omega)$,\ for\ $\alpha < \frac{N-4}{2}$.\
%
\end{lemma}
\emph{Proof}\, i)\ Let\ $u \in H_{0}^{2} (\Omega)$.\ Hardy's
inequality (\ref{rel1}) implies that
\[
\int_{\Omega} |x|^{-N} |v|^2\, dx = \int_{\Omega}
\frac{|u|^2}{|x|^4}\, dx < \infty.
\]
In this direction, from relations (\ref{3.1}), (\ref{3.6}) we
deduce that
\begin{eqnarray*}
\int_{\Omega} |x|^{-(N-4)} |\Delta v|^2\, dx \leq c_1 I[u] <
\infty,
\\ \int_{\Omega} |x|^{-(N-2)} |\nabla v|^2\, dx \leq c_2 I[u]<
\infty,
\end{eqnarray*}
for some positive constants\ $c_1,\; c_2$. Hence\ $||v||_{W}<
\infty$.

ii) Let\ $v \in C_{0}^{\infty} (\Omega)$\ and\ $u = |x|^{-\alpha}
v $,\ for\ $\alpha < \frac{N-4}{2}$.\ It is known (see
\cite{ft01}) that
\begin{equation} \label{6.2}
\biggl ( b - \frac{N-2}{2}  \biggr )^2 \int_{\Omega} |x|^{-2b-2}
w^2\, dx \leq \int_{\Omega} |x|^{-2b} |\nabla w|^2\,
\end{equation}
for any\ $w \in C^{\infty}_{0}$\ and\ $b \leq \frac{N-2}{2}$.\
Inequality (\ref{6.2}) for\ $b = a+1$,\ $a < \frac{N-4}{2}$,\
implies that
\begin{equation} \label{6.3}
\biggl ( a - \frac{N-4}{2}  \biggr )^2 \int_{\Omega} |x|^{-2a-4}
w^2\, dx \leq \int_{\Omega} |x|^{-2a-2} |\nabla w|^2.
\end{equation}
Hence, from (\ref{2.1}) and (\ref{6.3}) we conclude that
\begin{eqnarray*}
\int_{\Omega} |\Delta u|^2\, dx &\leq& C_{a} \biggl [
\int_{\Omega} |x|^{-2a} |\Delta v|^2\, dx + \int_{\Omega}
|x|^{-2a-4} |x \cdot \nabla v|^2\, dx +  \int_{\Omega}
|x|^{-2a-2} |\nabla v|^2\, dx \biggr ] \\
&\leq& C_{a} \biggl [ \int_{\Omega} |x|^{-(N-4)} |\Delta v|^2\, dx
+ \int_{\Omega} |x|^{-N} |x \cdot \nabla v|^2\, dx  +
\int_{\Omega} |x|^{-(N+2)} |\nabla v|^2\, dx \biggr ] \\
&\leq& C_{a}\; ||v||_{W^{2,2}_0 (\Omega, |x|^{-(N-4)})} < \infty
\end{eqnarray*}
and the proof is completed.\  $\blacksquare$
\begin{lemma} \label{l7.2}
The functionals\ $J,\; \mathbb{J}$\ are weakly lower
semicontinuous in\ $W^{2,2}_0 (\Omega, |x|^{-(N-4)})$.\
\end{lemma}
\emph{Proof}\, Let\ $v_n$\ be a weakly convergent sequence in\
$W^{2,2}_0 (\Omega, |x|^{-(N-4)})$,\ to some\ $v_0$.\ Assume also
the sequence\ $w_n := v_n -v_0$,\ with\ $w_n \rightharpoonup 0$.\
Then for\ $n$\ large enough, we may prove that
\begin{eqnarray*}
J(v_n) &=& J(w_n + v_0) = \\
&=& \int_{\Omega} |x|^{4-N} |\Delta w_n +\Delta v_0|^2\, dx -
N(N-4) \int_{\Omega} |x|^{-N} |x \nabla w_n + x \nabla v_0|^2\, dx
\\ &&+ \frac{N(N-4)}{2} \int_{\Omega} |x|^{2-N} |\nabla w_n +\nabla
v_0|^2\, dx = \\
&=& J(w_n) + J(v_0) + o (1).
\end{eqnarray*}
Since\ $J(w_n) \geq 0$,\ we conclude that\ $\liminf_{n \to \infty}
J(v_n) \geq J(v_0)$.\ The case of\ $\mathbb{J}$\ may be treated in
a similar way.\ $\blacksquare$ \vspace{0.2cm}
%
%
\subsection{Existence of minimizers for improved inequalities of (\ref{rel1})}
Assume the improved Hardy inequality
\begin{equation} \label{6.1.1}
\int_{\Omega} |\Delta u|^2\, dx \geq \biggl ( \frac{N(N-4)}{4}
\biggr )^2\, \int_{\Omega} \frac{u^2}{|x|^4}\, dx + b
\int_{\Omega} V\, u^2\, dx.
\end{equation}
We want the potential\ $V$\ to be a lower order potential compared
to the Hardy potential\ $\frac{1}{|x|^4}$.\ For that reason we
give the following definition of the admissible class\ ${\cal A}$\
of potentials :
%
We say that a potential\ $V$\ is an admissible potential, that is\
$V \in  {\cal A}$,\ if\ $V$\ is not everywhere nonpositive,\ $V
\in L^{\frac{N}{4}}_{loc}(\xO \setminus \{0\})$,\ and there exists
a positive constant\ $c$,\ such that
\begin{equation} \label{def1}
\int_{\xO} |\Delta u|^2 dx \geq  \left(\frac{N(N-4)}{4} \right)^2
    \int_{\xO} \frac{u^2}{|x|^4} dx +
 c \int_{\xO} |V|  u^2 dx,\;\;\; \mbox{for any}\;\; u \in H_{0}^2 (\Omega).
\end{equation}
%
The presence of the absolute value in the right hand side  of
(\ref{def1}) ensures that the negative part of $V$ is itself a
lower order potential compared to the Hardy potential, and
therefore the Hardy potential is truly present in (\ref{rel1}). As
a consequence of (\ref{4.11}), the class\ ${\cal A}$\ contains all
non everywhere nonpositive potentials\ $V$,\ such that\
$\int_{\xO} |V|^{\frac{N}{4}} X^{1-N/2} dx < \infty$.\
\vspace{0.2cm}

Actually, the best constants arising in the inequalities of type
(\ref{6.1.1}) in $u$-terms are equal with those ones arising in
the corresponding inequalities in $v$-terms
($v=|x|^{\frac{N-4}{2}} u$). For example, we have:
\begin{lemma} \label{l6.1.1}
The best constants
\begin{equation} \label{6.1.2}
c := \inf_{ \scriptsize \begin{array}{ll} & \;\;\; u \in H_{0}^{2} (\Omega), \vspace{0.1cm} \\
& \int_{\Omega} V\, u^2\, dx >0
\end{array}
  } \frac{I[u]}{\int_{\Omega} |V|\, u^2\, dx}
\end{equation}
and
\begin{equation} \label{6.1.3}
C := \inf_{ \scriptsize \begin{array}{ll} & \;\;\; v \in W^{2,2}_0 (\Omega, |x|^{-(N-4)}), \vspace{0.1cm} \\
& \int_{\Omega} |x|^{-(N-4)} V\, v^2\, dx >0
\end{array}
  } \frac{J(v)}{\int_{\Omega} |x|^{-(N-4)}  |V|\, v^2\, dx}
\end{equation}
are equal.
\end{lemma}
\emph{Proof}\, Let\ $c$,\ $C$\ be the best constants in
(\ref{6.1.2}) and (\ref{6.1.3}), respectively. For any\ $u \in
C_{0}^{\infty} (\Omega)$\ and\ $v = |x|^{\frac{N-4}{2}} u$,\ Lemma
\ref{l2.2} implies that
\[
\frac{I[u]}{\int_{\Omega} |V|\, u^2\, dx} =
\frac{J(v)}{\int_{\Omega} |x|^{-(N-4)}  |V|\, v^2\, dx}.
\]
Hence,\ $c \geq C$.\ Next we claim that\ $c \leq C$.\ Fix\
$\epsilon
>0$\ and assume the functions\ $v_{\epsilon} \in C_{0}^{\infty}
(\Omega)$,\ such that
\[
\frac{J(v_{\epsilon})}{\int_{\Omega} |x|^{-(N-4)}  |V|\,
v_{\epsilon}^2\, dx} \leq C + \epsilon.
\]
Let\ $0<a<\frac{N-4}{2}$.\ Lemma \ref{l6.1} implies that\ $u_{a,
\epsilon} = |x|^{-a} v_{\epsilon} \in H_{0}^{2} (\Omega)$\
providing that
\begin{equation}\label{6.1.4}
c \leq \frac{I[u_{a, \epsilon}]}{\int_{\Omega} |V|\, u_{a,
\epsilon}^2\, dx} = \frac{J_{a}(v_{\epsilon})}{\int_{\Omega}
|x|^{-2a-4}  |V|\, v_{\epsilon}^2\, dx},
\end{equation}
where
\begin{eqnarray*}
J_{a}(v) &:=& \int_{\Omega} |x|^{-2a} |\Delta v|^2\, dx - 4a(a+2)
\int_{\Omega} |x|^{-2a-4} |x \cdot \nabla v|^2\, dx \\
&& + 2a(a+2) \int_{\Omega} |x|^{-2a-2} |\nabla v|^2\, dx \\
&& + \biggl [ a(a+2)(-N+a+2)(-N+a+4) - \biggl ( \frac{N(N-4)}{4}
\biggr )^2 \biggr ] \int_{\Omega} |x|^{-2a-4} v^2\, dx.
\end{eqnarray*}
Next we calculate the limit of\ $J_{a}(v)$\ as\ $a \to
\frac{N-4}{2}^-$.\ It is clear that
\begin{eqnarray*}
\int_{\Omega} |x|^{-2a} |\Delta v|^2\, dx &\to& \int_{\Omega}
|x|^{-(N-4)} |\Delta v|^2\, dx, \\
\int_{\Omega} |x|^{-2a-2} |\nabla v|^2\, dx &\to& \int_{\Omega}
|x|^{-(N-2)} |\nabla v|^2\, dx.
\end{eqnarray*}
as\ $a \to \frac{N-4}{2}^-$.\ However, the problem arises in the
case of\ $\lim_{a \to \frac{N-4}{2}^-} \int_{\Omega} |x|^{-2a-4}
v^2\, dx$.\ In this case, we have that

\begin{eqnarray*}
\biggl [ a(a+2)(-N+a+2)(-N+a+4) - \biggl ( \frac{N(N-4)}{4} \biggr
)^2 \biggr ] \cdot \int_{\Omega} |x|^{-2a-4} v^2\, dx \\
\leq C\, ||v||_{\infty} \frac{a(a+2)(-N+a+2)(-N+a+4) - \biggl (
\frac{N(N-4)}{4} \biggr )^2}{N-2a-4} \to 0,
\end{eqnarray*}
as\ $a \to \frac{N-4}{2}^-$,\ hence
\[
\lim_{a \to \frac{N-4}{2}^-} J_{a}(v) = J(v),
\]
for any\ $v \in C_{0}^{\infty}(\Omega / \{0\})$.\
Taking the limit\ $a \to \frac{N-4}{2}^-$\ in (\ref{6.1.4}), we
obtain that
\[
c \leq C + \epsilon,
\]
for any fixed\ $\epsilon>0$\ and the proof is completed.\
$\blacksquare$ \vspace{0.2cm}

By the same argument the Hardy-Sobolev inequality (\ref{4.11})
takes the following form:
\begin{lemma} \label{l5.4}
Let\ $D \geq \sup_{x \in \xO} |x|$.\ Then, there exists\ $c>0$,\
such that
\begin{equation}
J[v] \geq c \left( \int_{\xO} |x|^{-N} |v|^{\frac{2N}{N-4}}
X^{\frac{2N-4}{N-2}}(\frac{|x|}{D}) dx \right)^{\frac{N-4}{N}},
\end{equation}
for every\ $v \in W^{2,2}_0 (\Omega, |x|^{-(N-4)})$.\
\end{lemma}

Define the following quantity
\[
Q[v] := \frac{J(v)}{\int_{\Omega} |x|^{-(N-4)}  V\, v^2\, dx}
\]
and set
\begin{equation} \label{6.1.5}
B := \inf_{ \scriptsize \begin{array}{ll} & \;\;\; v \in C_{0}^{\infty} (\Omega), \vspace{0.1cm} \\
& \int_{\Omega} |x|^{-(N-4)} V\, v^2\, dx >0
\end{array}
  } Q[v] = \inf_{ \scriptsize \begin{array}{ll} & \;\;\; v \in W^{2,2}_0 (\Omega, |x|^{-(N-4)}), \vspace{0.1cm}
  \\
& \int_{\Omega} |x|^{-(N-4)} V\, v^2\, dx >0
\end{array}
  } Q[v]
\end{equation}
By practically the same arguments as in Lemma \ref{l6.1.1} we have
that
\begin{lemma}
There holds:\ $B=b$
\end{lemma}

The local best constant of inequality (\ref{6.1.1}) can be written
as:
\begin{equation} \label{6.1.6}
C^0 := \lim_{r \downarrow 0} C_r,\;\;\;\;\;\;
C_r := \inf_{ \scriptsize \begin{array}{ll} & \;\;\; v \in C_{0}^{\infty} (B_r), \vspace{0.1cm} \\
& \int_{B_r} |x|^{-(N-4)} V\, v^2\, dx >0
\end{array}
  } \frac{J(v)}{\int_{B_r} |x|^{-(N-4)}  V\, v^2\, dx}.
\end{equation}
If there is no\ $v \in C_{0}^{\infty} (B_r)$,\ such that\
$\int_{B_r} |x|^{-(N-4)} V\, v^2\, dx >0$,\ for some\ $r>0$,\ we
set\ $C_r = \infty$.\ Observe that\ $B \leq C^0$.\
\setcounter{theorem}{6}
\begin{theorem} \label{p6.1.1}
Let
\begin{equation} \label{6.1.7}
B<C^0.
\end{equation}
Then\ $B$\ is achieved by some\ $v_0 \in W^{2,2}_0 (\Omega,
|x|^{-(N-4)})$.\
\end{theorem}
\emph{Proof}\, Let\ $\{ v_k \} \subset W^{2,2}_0 (\Omega,
|x|^{-(N-4)})$\ be a minimizing sequence for (\ref{6.1.5}),\ such
that
\begin{equation} \label{6.1.8}
\int_{\Omega} |x|^{-(N-4)}  V\, v_k^2\, dx =1,
\end{equation}
for every\ $k$.\ Hence\ $J(v_k) \to B$.\ Since\ $J(v_k)$\ is
bounded, from (\ref{3.1}) and (\ref{3.6}) we deduce that\ $\{ v_k
\}$\ must be bounded too, in\ $W^{2,2}_0 (\Omega, |x|^{-(N-4)})$.\
Therefore, there exists a subsequence, still denoted by\ $\{ v_k
\}$,\ such that
\[
v_k \rightharpoonup v_0,\;\;\; \mbox{in}\;\; W^{2,2}_0 (\Omega,
|x|^{-(N-4)})
\]
and
\[
v_k \to v_0,\;\;\; \mbox{in}\;\; L^2(\Omega/B_\rho),\;\;\;
\mbox{for every}\;\; \rho>0,
\]
for some\ $v_0 \in W^{2,2}_0 (\Omega, |x|^{-(N-4)})$.\ We set\
$w_k := v_k -v_0$.\ Then from (\ref{6.1.7}) we have that
\begin{equation}\label{6.1.9}
1=\int_{\Omega} |x|^{-(N-4)}  V\, w_k^2\, dx + \int_{\Omega}
|x|^{-(N-4)}  V\, v_0^2\, dx +o(1).
\end{equation}
In addition from Lemma \ref{l7.2} we deduce that
\[
B = J(w_k) + J(v_0) + o(1)
\]
or
\begin{equation}\label{6.1.10}
B \geq J(w_k) + B \int_{\Omega} |x|^{-(N-4)}  V\, v_0^2\, dx +
o(1)
\end{equation}
and
\begin{equation}\label{6.1.11}
B \geq J(v_0).
\end{equation}
Observe also that (\ref{6.1.7}) implies the existence of a\
$\rho>0$,\ sufficiently small, such that
\begin{equation}\label{6.1.12}
B \leq C_\rho = \inf_{ \scriptsize \begin{array}{ll} & \;\;\; v \in C_{0}^{\infty} (B_\rho), \vspace{0.1cm} \\
& \int_{B_\rho} |x|^{-(N-4)} V\, v^2\, dx >0
\end{array}
  } \frac{J(v)}{\int_{B_\rho} |x|^{-(N-4)}  V\, v^2\, dx}.
\end{equation}
Assume the cutoff function\ $\phi \in C_{0}^{\infty} (B_{\rho})$,\
such that\ $0 \leq \phi \leq 1$,\ in\ $B_{\rho}$\ and\ $\phi
\equiv 1$,\ in\ $B_{\rho/2}$.\ Set\ $w_k = \phi w_k + (1-\phi)
w_k$.\ Making some calculations we have that
\begin{eqnarray*}
J(w_k) &=& J(\phi w_k) + J((1-\phi)w_k) + 2 \int_{B_\rho}
|x|^{-(N-4)} \Delta (\phi w_k)\, \Delta
((1-\phi)w_k)\, dx \\
&& -N(N-4) \biggl [ 2 \int_{B_\rho} |x|^{-N} (x \cdot \nabla (\phi
w_k))\, (x \cdot \nabla ((1-\phi) w_k))\, dx \\
&&\;\;\;\;\;\;\;\;\;\;\;\;\;\;\;\;\;\;\;\;\; - \frac{1}{2}
\int_{B_\rho} |x|^{-N} \nabla (\phi w_k) \cdot \nabla ((1-\phi)
w_k)\, dx \biggr ] \\
&=& J(\phi w_k) + J((1-\phi)w_k) + 2 \int_{B_\rho} |x|^{-(N-4)}
\phi (1-\phi)\, |\Delta w_k|^2\, dx + o(1).
\end{eqnarray*}
Since\ $J((1-\phi)w_k) \geq 0$\ we obtain that
\begin{equation}\label{6.1.13}
J(w_k) \geq J(\phi w_k) + o(1).
\end{equation}
From (\ref{6.1.12}) we have that
\begin{equation} \label{6.1.14}
J(\phi w_k) \geq C_\rho \int_{B_\rho} |x|^{-(N-4)} V (\phi
w_k)^2\, dx.
\end{equation}
Since\ $V \in L^{N/4}_{loc} (\Omega/\{0\})$\ holds that
\begin{equation}\label{6.1.15}
\int_{\Omega/B_{\rho/2}} |x|^{-(N-4)} V w_k^2\, dx \to 0,\;\;\;
\mbox{as}\;\; k \to \infty.
\end{equation}
So, inequalities (\ref{6.1.13}), (\ref{6.1.14}) and (\ref{6.1.15})
imply that
\begin{equation}\label{6.1.16}
J(w_k) \geq C_\rho \int_{\Omega} |x|^{-(N-4)} V w_k^2\, dx + o(1).
\end{equation}
Then, from (\ref{6.1.9}) and (\ref{6.1.16}) we derive that
\[
J(w_k) \geq C_\rho \biggl ( 1 - \int_{\Omega} |x|^{-(N-4)} V
v_0^2\, dx \biggr ) + o(1).
\]
Taking into account (\ref{6.1.10}) we conclude that
\[
B \geq C_\rho \biggl ( 1 - \int_{\Omega} |x|^{-(N-4)} V v_0^2\, dx
\biggr ) + B \int_{\Omega} |x|^{-(N-4)}  V\, v_0^2\, dx + o(1),
\]
or
\[
(B-C_\rho) \biggl ( 1 - \int_{\Omega} |x|^{-(N-4)} V v_0^2\, dx
\biggr ) \geq 0
\]
which implies that
\[
\int_{\Omega} |x|^{-(N-4)} V v_0^2\, dx \geq 1
\]
and from (\ref{6.1.11}) that
\[
0 \leq \frac{J(v_0)}{\int_{\Omega} |x|^{-(N-4)} V v_0^2\, dx} \leq
B.
\]
It follows that\ $B$\ is attained by\ $v_0$.\ We note that
\[
\int_{\Omega} |x|^{-(N-4)} V v_0^2\, dx = 1
\]
and it follows from (\ref{6.1.10}) that\ $v_k$\ converges strongly
in\ $W^{2,2}_0 (\Omega, |x|^{-(N-4)})$\ to\ $v_0$.\ \
$\blacksquare$ \vspace{0.2cm}

We next look for an improvement of inequality (\ref{6.1.1}). That
is, for an inequality of the form:
\begin{equation} \label{5.22}
\int_{\xO} |\Delta u|^2 dx \geq \biggl ( \frac{N(N-4)}{4} \biggr
)^2 \int_{\xO} \frac{u^2}{|x|^4} dx + b \int_{\xO} V u^2 dx + b_1
\int_{\xO} W u^2 dx,\;\;\;\;\;\; u \in H_0^2 (\xO),
\end{equation}
where\ $V$\ and\ $W$\ are both in\ ${\cal A}$.\ Assuming that
(\ref{5.22}) holds true, the best constant\ $b_1$,\ is clearly
given by:
\begin{equation} \label{5.23}
b_1 = \inf_{\scriptsize
\begin{array}{c}
               u \in H_0^2 (\xO) \\
              \int_{\xO} W u^2 dx>0
               \end{array}}
\frac{ I[u] -b \int_{\xO} V  u^2 dx } {\int_{\xO} W u^2 dx}.
\end{equation}
By the same argument as in Lemma \ref{l6.1.1}, the constant\
$b_1$\ is also equal to:
\begin{equation} \label{5.24}
B_1 = \inf_{\scriptsize
\begin{array}{c}
               v \in W^{2,2}_0 (\Omega, |x|^{-(N-4)}) \\
              \int_{\xO} |x|^{-(N-4)}  W v^2 dx>0
               \end{array}}
\frac{J[v] - b \int_{\xO} |x|^{-(N-4)} V v^2 dx} {\int_{\xO}
|x|^{-(N-4)} W v^2 dx}.
\end{equation}
Notice that by the properties of $b=B$  we always have that $b_1
\geq 0$.
Conversely, if one defines\ $b_1 \geq 0$\ by (\ref{5.24}) it is
immediate that inequality (\ref{5.22}) holds true with\ $b_1$\
being the best constant. But of course, for (\ref{5.22}) to be an
improvement of the original inequality, we need\ $b_1$\ to be
strictly positive.
Our next result is a direct consequence of Proposition
\ref{p6.1.1} and provides conditions under which the original
inequality cannot be improved.
\setcounter{theorem}{7}
\begin{theorem} \label{p5.2} Suppose that\ $b < C^{0}$.\
Let\ $V$\ and\ $W$\ be both in\ ${\cal A}$.\ If\ $\phi$\ is the
minimizer of the quotient (\ref{6.1.5}) and
\[
\int_{\xO} |x|^{-(N-4)} W \phi^2 dx  >0,
\]
then\ $b_1=0$,\ that is, there is no further improvement of
(\ref{6.1.1}).
\end{theorem}
{\emph Proof}\, By our assumptions,\ $v=\phi$\ is an admissible
function in (\ref{5.24}). Moreover, for\ $v=\phi$\ the numerator
of (\ref{5.24}) becomes zero. In view of the fact that\ $b_1 \geq
0$,\ we conclude that\ $b_1 = 0$.\  $\blacksquare$ \vspace{0.2cm}

It follows in particular that if $W \geq 0$,  we cannot improve
(\ref{rel1}). Thus, the following result has been proved.
\setcounter{theorem}{8}
\begin{theorem} \label{t5.6}
Let\ $V \in {\cal A}$.\ If
\[
b < C^{0},
\]
then, we cannot improve (\ref{6.1.1}) by adding a nonnegative
potential\ $W \in {\cal A}$.\
\end{theorem}

As a consequence of (\ref{4.11}) and Theorem \ref{t5.6} we have:
\begin{coro} \label{c5.7}
Let\ $D >  \sup_{x \in \xO} |x|$.\ Suppose\ $V$\ is not everywhere
nonpositive and such that
\[  \int_{\xO} |V|^{\frac{N}{4}}
X^{1-N/2}(|x|/D)\; dx < \infty.
\]
Then,\ $V \in {\cal A}$\ but there is no further improvement of
(\ref{6.1.1}) with a nonnegative $W \in {\cal A} $.
\end{coro}
{\emph Proof}\, Applying Holder's inequality we get:
\[
\int_{\xO} |x|^{-(N-4)} |V|  v^2\; dx \leq \left( \int_{\xO}
|V|^{\frac{N}{4}} X^{1-N/2}\; dx \right)^{\frac{4}{N}}
\left(\int_{\xO} |x|^{-N} X^{\frac{2N-4}{N-4}}
|v|^{\frac{2N}{N-4}}\; dx \right)^{\frac{N-4}{N}}.
\]
The first integral is bounded by our assumption, whereas the
second integral is bounded  from above  by\ $ C\; J[v]$\ (cf Lemma
\ref{l5.4}). Thus we proved that\ $V \in {\cal A}$.\ Using once
more Holder's inequality in\ $B_r$\ and the definition of\ $C_r$\
(cf (\ref{6.1.6})) we easily see that:
\[
C_r  \geq \frac{C}{\left(\int_{B_r} |V|^{\frac{N}{4}} X^{1-N/2}\;
dx \right)^{\frac{4}{N}}} \to \infty,\;\;\;\;\;\; {\rm as}\;\;\; r
\to 0,
\]
whence\ ${\cal  C}^{0} = +\infty$.\ Thus, all conditions of
Theorem \ref{t5.6} are satisfied and the result follows.\
$\blacksquare$ \vspace{0.2cm}
%
%
\subsection{Existence of minimizers for improved inequalities of (\ref{rel2})}
Assume the following improved inequality
\begin{equation}\label{6.2.1}
\int_{\Omega} |\Delta u|^2\, dx \geq \frac{N^2}{4}\, \int_{\Omega}
\frac{|\nabla u|^2}{|x|^2}\, dx + \mathsf{b} \int_{\Omega} V\,
|\nabla u|^2\, dx.
\end{equation}
We want the potential\ $V$\ to be a lower order potential compared
to the Hardy potential\ $\frac{1}{|x|^2}$.\ For that reason we
give the following definition of the admissible class\
${\mathbb{A}}$\ of potentials :
\begin{defin}
We say that a potential\ $V$\ is an admissible potential, that is\
$V \in  {\mathbb{A}}$,\ if\ $V$\ is not everywhere nonpositive,\
$V \in L^{\frac{N}{2}}_{loc}(\xO \setminus \{0\})$,\ and there
exists a positive constant\ $c$,\ such that
\begin{equation} \label{def2}
\int_{\xO} |\Delta u|^2 dx \geq  \frac{N^2}{4} \int_{\xO}
\frac{|\nabla u|^2}{|x|^2} dx + c \int_{\xO} |V|  |\nabla u|^2
dx,\;\;\; \mbox{for any}\;\; u \in H_{0}^2 (\Omega).
\end{equation}
\end{defin}
The presence of the absolute value in the right hand side  of
(\ref{def1}) ensures that the negative part of $V$ is itself a
lower order potential compared to the Hardy potential, and
therefore the Hardy potential is truly present in (\ref{rel2}). As
a consequence of (\ref{4.20}), the class\ $\mathbb{A}$\ contains
all non everywhere nonpositive potentials\ $V$,\ such
that\ $\int_{\xO} |V|^{\frac{N}{2}} X^{1-N} dx < \infty$.\ 
\begin{lemma} \label{l6.2.1}
The best constants
\begin{equation} \label{6.2.2}
\mathsf{c} := \inf_{ \scriptsize \begin{array}{ll} & \;\;\; u \in H_{0}^{2} (\Omega), \vspace{0.1cm} \\
& \int_{\Omega} V\, |\nabla u|^2\, dx >0
\end{array}
  } \frac{\mathbb{I}[u]}{\int_{\Omega} |V|\, |\nabla u|^2\, dx},
\end{equation}
and
\begin{equation} \label{6.2.3}
\mathbb{C} := \inf_{ \scriptsize \begin{array}{ll} & \;\;\; v \in W^{2,2}_0 (\Omega, |x|^{-(N-4)}), \vspace{0.1cm} \\
& \int_{\Omega} |x|^{-(N-4)} V\, |\nabla v - \frac{N-4}{2}
\frac{x}{|x|^2} v|^2\, dx
>0
\end{array}
  } \frac{\mathbb{J}(v)}{\int_{\Omega} |x|^{-(N-4)}  |V|\, | \nabla v - \frac{N-4}{2} \frac{x}{|x|^2} v |^2\,
  dx},
\end{equation}
are equal.
\end{lemma}
\emph{Proof}\, Let\ $\mathsf{c}$,\ $\mathbb{C}$\ be the best
constants in (\ref{6.2.2}) and (\ref{6.2.3}), respectively. For
any\ $u \in C_{0}^{\infty} (\Omega)$\ and\ $v =
|x|^{\frac{N-4}{2}} u$,\ Lemma \ref{l2.2} implies that
\[
\frac{\mathbb{I}[u]}{\int_{\Omega} |V|\, |\nabla u|^2\, dx} =
\frac{\mathbb{J}(v)}{\int_{\Omega} |x|^{-(N-4)}  |V|\,
\tilde{v}^2\, dx},
\]
hence\ $\mathsf{c} \geq \mathbb{C}$.\ Next we claim that\
$\mathsf{c} \leq \mathbb{C}$.\ Fix\ $\epsilon >0$\ and assume the
functions\ $v_{\epsilon} \in C_{0}^{\infty} (\Omega)$,\ such that
\[
\frac{\mathbb{J}(v_{\epsilon})}{\int_{\Omega} |x|^{-(N-4)}  |V|\,
|\nabla v_{\epsilon} - \frac{N-4}{2} \frac{x}{|x|^2} v_{\epsilon}
|^2\, dx} \leq \mathbb{C} + \epsilon.
\]
Let\ $0<a<\frac{N-4}{2}$.\ Lemma \ref{l6.1} implies that\ $u_{a,
\epsilon} = |x|^{-a} v_{\epsilon} \in H_{0}^{2} (\Omega)$\
providing that
\begin{equation}\label{6.2.4}
\mathbb{C}_1 \leq \frac{\mathbb{I}[u_{a, \epsilon}]}{\int_{\Omega}
|V|\, | \nabla u_{a, \epsilon}|^2\, dx} =
\frac{\mathbb{J}_{a}(v_{\epsilon})}{\int_{\Omega} |x|^{-2a} |V|\,
|\nabla v_{\epsilon} - a \frac{x}{|x|^2} v_{\epsilon}|^2\, dx},
\end{equation}
where
\begin{eqnarray*}
\mathbb{J}_{a}(v) &:=& \int_{\Omega} |x|^{-2a} |\Delta v|^2\, dx -
4a(a+2)
\int_{\Omega} |x|^{-2a-4} |x \cdot \nabla v|^2\, dx \\
&& + \biggl [ 2a(a+2) - \frac{N^2}{4} \biggr ] \int_{\Omega}
|x|^{-2a-2} |\nabla v|^2\, dx \\ && + a(-N+a+4) \biggl [
(a+2)(-N+a+2) + \frac{N^2}{4} \biggr ] \int_{\Omega} |x|^{-2a-4}
v^2\, dx.
\end{eqnarray*}
Following similar arguments as in Lemma \ref{l6.1.1} we may prove
that
\[
a(-N+a+4) \biggl [ (a+2)(-N+a+2) + \frac{N^2}{4} \biggr ] \cdot
\int_{\Omega} |x|^{-2a-4} v^2\, dx \to 0,
\]
as\ $a \to \frac{N-4}{2}^{-}$,\ hence
\begin{equation} \label{6.4}
\lim_{a \to \frac{N-4}{2}^{-}} \mathbb{J}_{a}(v) = \mathbb{J}(v),
\end{equation}
for any\ $v \in C_{0}^{\infty}(\Omega / \{0\})$.\ Using now
(\ref{6.4}) and
\[
\lim_{a \to \frac{N-4}{2}^{-}} \int_{\Omega} |x|^{-2a} |V|\,
|\nabla v_{\epsilon} - a \frac{x}{|x|^2} v_{\epsilon}|^2\, dx =
\int_{\Omega} |x|^{-(N-4)} |V|\, |\nabla v_{\epsilon} -
\frac{N-4}{2} \frac{x}{|x|^2} v_{\epsilon}|^2\, dx
\]
we obtain that
\[
\mathsf{c} \leq \mathbb{C} + \epsilon,
\]
for any fixed\ $\epsilon>0$\ and the proof is completed.\
$\blacksquare$ \vspace{0.2cm}

By the same argument the Hardy-Sobolev inequality (\ref{4.20})
takes the following form:
\begin{lemma} \label{l5.4b}
Let\ $D \geq \sup_{x \in \xO} |x|$.\ Then, there exists\ $c>0$,\
such that
\begin{equation}
\mathbb{J}[v] \geq c \left( \int_{\xO} |x|^{-\frac{N(N-4)}{N-2}}
|\nabla v - \frac{N-4}{2} \frac{x}{|x|^2} v|^{\frac{2N}{N-2}}
X^{\frac{2N-2}{N-2}}(\frac{|x|}{D}) dx \right)^{\frac{N-2}{N}},
\end{equation}
for every\ $v \in W^{2,2}_0 (\Omega, |x|^{-(N-4)})$.\
\end{lemma}

Define the following quantity
\[
\mathbb{Q}[v] := \frac{\mathbb{J}(v)}{\int_{\Omega} |x|^{-(N-4)}
V\, |\nabla v - \frac{N-4}{2} \frac{x}{|x|^2} v|^2\, dx}
\]
and set
\begin{equation} \label{6.2.5}
\mathbb{B} := 
\inf_{ \scriptsize \begin{array}{ll} & \;\;\; v \in W^{2,2}_0
(\Omega, |x|^{-(N-4)}), \vspace{0.1cm}
  \\
& \int_{\Omega} |x|^{-(N-4)} V\, |\nabla v - \frac{N-4}{2}
\frac{x}{|x|^2} v|^2\, dx >0
\end{array}
  } \mathbb{Q}[v].
\end{equation}
%

By practically the same arguments as in Lemma \ref{l6.1.1} we have
that
\begin{lemma}
There holds:\ $\mathbb{B}=b$
\end{lemma}

The local best constant of inequality (\ref{6.2.1}) can be written
as:
\begin{equation} \label{6.2.6}
\mathbb{C}^0 := \lim_{r \downarrow 0} \mathbb{C}_r,\;\;
\mathbb{C}_r := \inf_{ \scriptsize \begin{array}{ll} &
\;\;\;\;\;\;\;\;\;\;\;\;\;\;\;\;\;\;\;\;\;\;\;
v \in C_{0}^{\infty} (B_r), \vspace{0.1cm} \\
& \int_{B_r} |x|^{-(N-4)} V |\nabla v - \frac{N-4}{2}
\frac{x}{|x|^2} v|^2 dx >0
\end{array}
  } \frac{\mathbb{J}(v)}{\int_{B_r} |x|^{-(N-4)}  V\, |\nabla v -
  \frac{N-4}{2} \frac{x}{|x|^2} v|^2\, dx}.
\end{equation}
If there is no\ $v \in C_{0}^{\infty} (B_r)$,\ such that\
$\int_{B_r} |x|^{-(N-4)} V\, |\nabla v -
  \frac{N-4}{2} \frac{x}{|x|^2} v|^2\, dx >0$,\ for some\
$r>0$,\ we set\ $\mathbb{C}_r = \infty$.\ Observe that\
$\mathbb{B} \leq \mathbb{C}^0$.\

We introduce the space\ $\mathcal{V}$\ to be defined as the
closure of the\ $C_{0}^{\infty}(\Omega)$\ functions with respect
to the norm
\begin{eqnarray} \label{6.2.7}
||u||^{2}_{\mathcal{V}} &:=& \int_{\Omega} |x|^{-(N-4)} |\Delta
v|^2\, dx + \int_{\Omega} |x|^{-N} |x \cdot \nabla v|^2\, dx +
\int_{\Omega} |x|^{-(N-2)} |\nabla v|^2\, dx + \int_{\Omega}
|x|^{-N}\, v^2\, dx \nonumber \\ && + \int_{\Omega} |x|^{-(N-4)}\,
|V|\, \biggl | \nabla v -\frac{N-4}{2} \frac{x}{|x|^2} v \biggr
|^2\, dx.
\end{eqnarray}
It is clear that\ $\mathcal{V}$\ is a Hilbert space with inner
product
\begin{eqnarray*}
< \phi, \psi >_{\mathcal{V}} &:=& \int_{\Omega} |x|^{-(N-4)}
\Delta \phi\, \Delta \psi\, dx + \int_{\Omega} |x|^{-N} (x \cdot
\nabla \phi) (x \cdot \nabla \psi)\, dx + \int_{\Omega}
|x|^{-(N-2)} \nabla \phi \cdot \nabla \psi\, dx \nonumber
\\ && + \int_{\Omega}
|x|^{-N} \phi\, \psi\, dx + \int_{\Omega} |x|^{-(N-4)} |V| (\nabla
\phi -\frac{N-4}{2} \frac{x}{|x|^2} \phi)(\nabla \psi
-\frac{N-4}{2} \frac{x}{|x|^2} \psi),
\end{eqnarray*}
satisfying
\begin{equation} \label{6.2.8}
||v||_{W} \leq ||v||_{\mathcal{V}} \leq c_0 ||v||_{W},\;\;\;
\mbox{for any}\;\; v \in \mathcal{V}.
\end{equation}
\setcounter{theorem}{14}
\begin{theorem} \label{p6.2.1}
Let
\begin{equation} \label{6.2.9}
\mathbb{B}<\mathbb{C}^0.
\end{equation}
Then\ $\mathbb{B}$\ is achieved by some\ $v_0 \in W^{2,2}_0
(\Omega, |x|^{-(N-4)})$.\
\end{theorem}
\emph{Proof}\, Let\ $\{ v_k \} \subset W^{2,2}_0 (\Omega,
|x|^{-(N-4)})$\ be a minimizing sequence for (\ref{6.2.1}),\ such
that
\begin{equation} \label{6.2.10}
L(v_k) := \int_{\Omega} |x|^{-(N-4)}  V\, \biggl | \nabla v_k
-\frac{N-4}{2} \frac{x}{|x|^2} v_k \biggr |^2\, dx =1,
\end{equation}
for every\ $k$.\ Hence\ $\mathbb{J}(v_k) \to B$.\ Since\
$\mathbb{J}(v_k)$\ is bounded, from (\ref{3.6}) and (\ref{3.7}) we
deduce that\ $\{ v_k \}$\ must be bounded too, in\ $W^{2,2}_0
(\Omega, |x|^{-(N-4)})$.\ Therefore, there exists a subsequence,
still denoted by\ $\{ v_k \}$,\ such that
\[
v_k \rightharpoonup v_0,\;\;\; \mbox{in}\;\; W^{2,2}_0 (\Omega,
|x|^{-(N-4)})
\]
and
\[
v_k \to v_0,\;\;\; \mbox{in}\;\; L^2(\Omega/B_\rho),\;\;\;
\mbox{for every}\;\; \rho>0,
\]
for some\ $v_0 \in W^{2,2}_0 (\Omega, |x|^{-(N-4)})$.\ We set\
$w_k := v_k -v_0$.\ Then from (\ref{6.2.8}) and (\ref{6.2.10}) we
have that
\begin{equation}\label{6.2.11}
1= L(w_k + v_0) = L(w_k) + L(v_0) +o(1).
\end{equation}
Following now the same steps as in the proof of Proposition
\ref{p6.1.1} we conclude that\ $L(v_0) \geq 1$,\ hence
\[
\frac{\mathbb{J}}{L(v_0)} \geq \mathbb{B}.
\]
This last inequality implies that\ $\mathbb{B}$\ is attained by\
$v_0$,\ such that\ $L(v_0) =1$\ and the proof is completed.\
$\blacksquare$ \vspace{0.2cm}

We next look for an improvement of inequality (\ref{6.2.1}). That
is, for an inequality of the form:
\begin{equation} \label{5.22b}
\int_{\xO} |\Delta u|^2 dx \geq \frac{N^2}{4} \int_{\xO}
\frac{|\nabla u|^2}{|x|^4} dx + \mathsf{b} \int_{\xO} V |\nabla
u|^2 dx + \mathsf{b}_1 \int_{\xO} W |\nabla u|^2 dx,\;\;\;\;\;\; u
\in H_0^2 (\xO),
\end{equation}
where\ $V$\ and\ $W$\ are both in\ $\mathbb{A}$.\ Assuming that
(\ref{5.22b}) holds true, the best constant\ $b_1$,\ is clearly
given by:
\begin{equation} \label{5.23b}
\mathsf{b}_1 = \inf_{\scriptsize
\begin{array}{c}
               u \in H_0^2 (\xO) \\
              \int_{\xO} W |\nabla u|^2 dx>0
               \end{array}}
\frac{ \mathbb{I}[u] -b \int_{\xO} V  |\nabla u|^2 dx }
{\int_{\xO} W |\nabla u|^2 dx}.
\end{equation}
By the same argument as in Lemma \ref{l6.2.1}, the constant\
$\mathsf{b}_1$\ is also equal to:
\begin{equation} \label{5.24b}
\mathbb{B}_1 = \inf_{\scriptsize
\begin{array}{c}
               v \in W^{2,2}_0 (\Omega, |x|^{-(N-4)}) \\
              \int_{\xO} |x|^{-(N-4)}  W |\nabla v -
  \frac{N-4}{2} \frac{x}{|x|^2} v|^2 dx>0
               \end{array}}
\frac{\mathbb{J}[v] - b \int_{\xO} |x|^{-(N-4)} V |\nabla v -
  \frac{N-4}{2} \frac{x}{|x|^2} v|^2 dx}
{\int_{\xO} |x|^{-(N-4)} W |\nabla v -
  \frac{N-4}{2} \frac{x}{|x|^2} v|^2 dx}.
\end{equation}
Notice that by the properties of $\mathsf{b}=\mathbb{B}$  we
always have that $\mathsf{b}_1 \geq 0$.
Conversely, if one defines\ $\mathsf{b}_1 \geq 0$\ by
(\ref{5.24b}) it is immediate that inequality (\ref{5.22b}) holds
true with\ $\mathsf{b}_1$\ being the best constant. But of course,
for (\ref{5.22b}) to be an improvement of the original inequality,
we need\ $\mathsf{b}_1$\ to be strictly positive.
Our next result is a direct consequence of Proposition
\ref{p6.2.1} and provides conditions under which the original
inequality cannot be improved.
\setcounter{theorem}{15}
\begin{theorem} \label{p5.2b} Suppose that\
$\mathbb{B}<\mathbb{C}^0$.\ Let\ $V$\ and\ $W$\ be both in\
$\mathbb{A}$.\ If\ $\phi$\ is the minimizer of the quotient
(\ref{6.2.5}) and
\[
\int_{\xO} |x|^{-(N-4)} W |\nabla \phi -
  \frac{N-4}{2} \frac{x}{|x|^2} \phi|^2 dx  >0,
\]
then\ $\mathsf{b}_1=0$,\ that is, there is no further improvement
of (\ref{6.2.1}).
\end{theorem}
{\emph Proof}\, By our assumptions,\ $v=\phi$\ is an admissible
function in (\ref{5.24b}). Moreover, for\ $v=\phi$\ the numerator
of (\ref{5.24b}) becomes zero. In view of the fact that\
$\mathsf{b}_1 \geq 0$,\ we conclude that\ $\mathsf{b}_1 = 0$.\
$\blacksquare$ \vspace{0.2cm}

It follows in particular that if $W \geq 0$,  we cannot improve
(\ref{6.2.1}). Thus, the following result has been proved.
\setcounter{theorem}{16}
\begin{theorem} \label{t5.6b}
Let\ $V \in \mathbb{A}$.\ If
\[
\mathbb{B}<\mathbb{C}^0,
\]
then, we cannot improve (\ref{6.2.1}) by adding a nonnegative
potential\ $W \in \mathbb{A}$.\
\end{theorem}

As a consequence of (\ref{4.20}) and Theorem \ref{t5.6b} we have:
\begin{coro} \label{c5.7b}
Let\ $D >  \sup_{x \in \xO} |x|$.\ Suppose\ $V$\ is not everywhere
nonpositive and such that
\[
\int_{\xO} |V|^{\frac{N}{2}} X^{1-N}(|x|/D)\; dx < \infty.
\]
Then,\ $V \in \mathbb{A}$\ but there is no further improvement of
(\ref{6.2.1}) with a nonnegative $W \in \mathbb{A} $.
\end{coro}
{\emph Proof}\, Applying Holder's inequality we get:
\begin{eqnarray*}
\int_{\xO} |x|^{-(N-4)} |V|  |\nabla v - \frac{N-4}{2}
\frac{x}{|x|^2} v|^2\; dx \leq \left( \int_{\xO} |V|^{\frac{N}{2}}
X^{1-N}\; dx \right)^{\frac{2}{N}} \cdot\;\;\;\;\;\;\;\;\;\;\;\;\;
\;\;\;\;\;\;\;\;\;\;\;\;\;\;\;\;\;\;\;\;\;\;\;\;\;\;\; \\
\;\;\;\;\;\;\;\;\;\;\;\;\;\;\;\;\;\;\;\;\;\;\;\;\;\;\;\;\;\;
\;\;\;\;\;\;\;\;\;\;\;\;\;\;\;\;\;\;\;\;\;\;\;\;\;\;\; \cdot
\left(\int_{\xO} |x|^{-\frac{N(N-4)}{N-2}} X^{\frac{2N-2}{N-2}}
|\nabla v - \frac{N-4}{2} \frac{x}{|x|^2} v|^{\frac{2N}{N-2}}\; dx
\right)^{\frac{N-2}{N}}.
\end{eqnarray*}
The first integral is bounded by our assumption, whereas the
second integral is bounded  from above  by\ $ C\; \mathbb{J}[v]$\
(cf Lemma \ref{l5.4b}). Thus we proved that\ $V \in \mathbb{A}$.\
Using, as in the proof of Corollary \ref{c5.7}, Holder's
inequality in\ $B_r$\ and the definition of\ $\mathbb{C}_r$\ (cf
(\ref{6.2.6})) we easily get that\ $\mathbb{C}^{0} = +\infty$.\
Thus, all conditions of Theorem \ref{t5.6b} are satisfied and the
result follows.\ $\blacksquare$ \vspace{0.6cm}
%
%
%
%
\section{PART II. THE POLYHARMONIC OPERATOR}
\setcounter{equation}{0}
In this part we prove some improved Hardy-Rellich inequalities
involving the polyharmonic operator. More precisely, we give the
proof of the Theorems \ref{t8.4} to \ref{tp} for which, we have to
establish certain inequalities concerning (\ref{rel1m}) and
(\ref{rel2m}).
%
%
\subsection{The Inequality (\ref{rel1m})}
\setcounter{equation}{0}
\begin{lemma} \label{l8.1}
Suppose\ $N \geq 5$\ and\ $0 \leq m<\frac{N-4}{2}$.\ For any\ $u
\in C_{0}^{\infty} (\Omega)$,\ we set\ $v = |x|^{a} u$.\ Then, the
following equality holds.
\begin{eqnarray} \label{8.1}
\int_{\Omega} \frac{|\Delta u|^2}{|x|^{2m}}\, dx = \int_{\Omega}
|x|^{-2m-2a} |\Delta v|^2\, dx - 4a(2m+2+a) \int_{\Omega}
|x|^{-2a-4-2m} (x \cdot \nabla v)^2\, dx \nonumber \\ + 2a(a+2+2m)
\int_{\Omega} |x|^{-2a-2-2m} |\nabla v|^2\, dx \nonumber \\ +
\biggl(a^2(a+2-N)^2 -2a(a+2-N)(m+1)(N-4-2m-2a)\biggr)
\int_{\Omega} |x|^{-2a-4-2m} v^2\, dx,
\end{eqnarray}
\end{lemma}
\begin{lemma} \label{l8.2}
Suppose\ $N \geq 5$\ and\ $0 \leq m<\frac{N-4}{2}$.\ For any\ $u
\in C_{0}^{\infty} (\Omega)$,\ we set\ $v = |x|^{\frac{N-4-2m}{2}}
u$.\ Then, the following equalities hold.
\begin{eqnarray*}
i)\;\;\; \int_{\Omega} \frac{|\nabla u|^2}{|x|^{2m+2}}\, dx =
\int_{\Omega} |x|^{-(N-2)} |\nabla v|^2\, dx + \biggl (
\frac{N-4-2m}{2} \biggr )^2 \int_{\Omega} |x|^{-N} |v|^2\,
dx,\;\;\;\;\;\;\;\;\;\;\;\;\;\;\;\;\;\;\;\;
\end{eqnarray*}
\begin{eqnarray*}
ii)\;\;\; \int_{\Omega} \frac{|\Delta u|^2}{|x|^{2m}}\, dx -
\biggl ( \frac{(N+2m)(N-4-2m)}{4} \biggr )^2\, \int_{\Omega}
\frac{u^2}{|x|^{2m+4}}\, dx = \;\;\;\;\;\;\;\;\;\;\;\;\;\;\;
\;\;\;\;\;\;\;\;\;\;\;\;\;\;\;\;\;\;\;\;\;\;\;\;\; \\
\;\;\;\;\;\; \int_{\Omega} |x|^{-(N-4)} |\Delta v|^2\, dx -
(N+2m)(N-4-2m)
\int_{\Omega} |x|^{-N} (x \cdot \nabla v)^2\, dx \\
+ \frac{(N+2m)(N-4-2m)}{2} \int_{\Omega} |x|^{-(N-2)} |\nabla
v|^2\, dx,
\end{eqnarray*}
\end{lemma}

\setcounter{theorem}{2}
\begin{theorem} \label{p8.3}
Suppose\ $N \geq 5$\ and\ $0 \leq m<\frac{N-4}{2}$.\ For any\ $u
\in C_{0}^{\infty} (\Omega)$,\ we set\ $v = |x|^{\frac{N-4-2m}{2}}
u$.\ Then, the following inequality holds.
\begin{equation} \label{8.2}
\int_{\Omega} \frac{|\Delta u|^2}{|x|^{2m}}\, dx - \biggl (
\frac{(N+2m)(N-4-2m)}{4} \biggr )^2\, \int_{\Omega}
\frac{u^2}{|x|^{2m+4}}\, dx \geq A(N,m)\; \int_{\Omega}
|x|^{-(N-2)} |\nabla v|^2\, dx,
\end{equation}
where
\[
A(N,m):= \biggl\{
\begin{array}{ll}
(N-1) + \frac{1}{2} (N+2m)(N-4-2m),   &  m > \frac{-2+\sqrt{N-1}}{2}, \\
4(1+m)^2 + \frac{(N+2m)(N-4-2m)}{2} ,  & m \leq
\frac{-2+\sqrt{N-1}}{2},
\end{array}
\]
Moreover, the constant $4(1+m)^2 + \frac{(N+2m)(N-4-2m)}{2}$ for\
$m < \frac{-2+\sqrt{N-1}}{2}$ is the best.
\end{theorem}
\emph{Proof}\, We proceed by using Lemma \ref{l8.2} (ii) and
decomposing\ $v$\ into spherical harmonics. The equalities
(\ref{2.8})-(\ref{2.10}) imply that is enough to prove that
\begin{eqnarray}
\biggl [ (k+1)^2 + (2k+N-1)(N-3) - \frac{1}{2} (N+2m)(N-4-2m) - A
\biggr ] \int_{\mathbb{R}^N} r^{2k-N+2} |\nabla g_k|^2\, dx \nonumber \\
\label{8.3} \geq  \biggl [ Ak(N-2) - \frac{k}{2}
(N+2m)(N-4-2m)(2k+N-2) \biggr ] \int_{\mathbb{R}^N} r^{2k-N}
(g_k)^2\, dx.
\end{eqnarray}
For\ $k=0$\ from (\ref{8.3}) we obtain that
\[
A \leq A_0 \equiv (N-2)^2 - \frac{1}{2} (N+2m)(N-4-2m),
\]
while for\ $k \ne 0$\ we obtain that
\[
A \leq A_1 \equiv (N-1) + \frac{1}{2} (N+2m)(N-4-2m),
\]
which corresponds for\ $k=1$.\ Then, we conclude that\ $A$\ must
be the minimum of\ $A_0$,\ $A_1$,\ or\ $A=A(N,m)$.\ Let\ $m <
\frac{-2+\sqrt{N-1}}{2}$\ and consider the minimizing sequences
$u^{\epsilon}$\ and\ $v^{\epsilon}$\
\[
u^{\epsilon} := r^{-\frac{N-4}{2}+m+\epsilon}
X_{1}^{\frac{-1+a_1}{2}} \phi(r),\;\;\;\;\;\; v^{\epsilon} :=
r^{\frac{N-4}{2}-m} u^{\epsilon} = r^{\epsilon}
X_{1}^{\frac{-1+a_1}{2}} \phi(r),
\]
in a similar way as in Section 4. Then, we have that
\[
\frac{\int_{\Omega} \frac{|\Delta u^{\epsilon}|^2}{|x|^{2m}}\, dx
- \biggl ( \frac{(N+2m)(N-4-2m)}{4} \biggr )^2\, \int_{\Omega}
\frac{(u^{\epsilon})^2}{|x|^{2m+4}}\, dx}{\int_{\Omega}
|x|^{-(N-2)} |\nabla v^{\epsilon}|^2\, dx} \to 4(1+m)^2 +
\frac{(N+2m)(N-4-2m)}{2},
\]
as\ $\epsilon \to 0^+$\ and\ $a_1 \to 0^+$.\ $\blacksquare$
\vspace{0.2cm}

Observe that\ for\ $m=0$\ Proposition \ref{8.3} implies that
\[
A(N,0) =  \biggl ( 4 + \frac{N(N-4)}{2} \biggr ),
\]
which is the result stated in Proposition \ref{p3.1}.\
\vspace{0.3cm} \\
\emph{Proof of Theorem \ref{t8.4}}\, When\ $0<m \leq
\frac{-2+\sqrt{N-1}}{2}$,\ inequality (\ref{8.4}) is an immediate
consequence from Proposition \ref{p8.3} and Inequality
(\ref{1.3})\ However, we will establish this for the whole range
of\ $m \in [0, \frac{N-4}{4})$.\ Once more we do the change of
variable of (\ref{2.8})-(\ref{2.10}). Then the inequality will be
true provided we will establish the following inequality
\begin{eqnarray*}
\biggl [ (k+1)^2 + (2k+N-1)(N-3) - \frac{1}{2} (N+2m)(N-4-2m)
\biggr ] \int_{\mathbb{R}^N} r^{2k-N+2} |\nabla g_k|^2\, dx \nonumber \\
+ \frac{k}{2} (N+2m)(N-4-2m)(2k+N-2) \int_{\mathbb{R}^N} r^{2k-N}
(g_k)^2\, \sum_{i=1}^{k-1} X_1^2
 \ldots X_k^2\, dx \\
\geq  \biggl [ (1+m)^2 + \frac{(N+2m)(N-4-2m)}{8} \biggr ]
\int_{\mathbb{R}^N} r^{2k-N} (g_k)^2 \sum_{i=1}^{k-1} X_1^2
 \ldots X_k^2\, dx.
\end{eqnarray*}
However, the worst case is for\ $k=0$,\ but this follows from
\ref{1.5}.  To establish the best constants we will treat
initially the case $m \leq \frac{-2+\sqrt{N-1}}{2}$. The proof of
it follows the same lines as in Section 4. For this we fix small
parameters\ $\epsilon,\, a_1,\, a_2,\,...,a_k >0$\ and define
\[
u(x) := w(x)\, \phi(|x|),\;\;\;
w(x) := |x|^{-\frac{N-4}{2}+m+\epsilon} X_{1}^{\frac{-1+a_1}{2}}
X_{2}^{\frac{-1+a_2}{2}} \cdots X_{k}^{\frac{-1+a_k}{2}},
\]
where\ $X_l = X_1 (X_{l-1})$,\ $l=2,...,k$.\ and\ $\phi(r) \in
C^{\infty}_{0} (B_1)$\ is a smooth cutoff function, such that\ $0
\leq \phi \leq 1$, with\ $\phi \equiv 1$\ in\ $B_{1/2}$.\
Following similar arguments as in Section 4 we may obtain that
\[
\begin{array}{ccll}
\int_{\Omega} \frac{|\Delta u|^2}{|x|^{2m}}\, dx - \biggl (
\frac{(N+2m)(N-4-2m)}{4} \biggr )^2 \int_{\Omega}
\frac{u^2}{|x|^{2m+4}} dx-\;\;\;\;\;\;\;\;\;\;\;\;\;\;\;
\;\;\;\;\;\;\;\;\;\;\;\;\;\;\;\;\;\;\;\;\;\;\;\;\;\;\;\;\;\;
\;\;\;\;\;\;\;\;\;\;\;\;\;\;\;\;\;\;\;\;\;\;\;\;\;\;\;\;\;\;
\\
\;\;\;\;\;\; \biggl ( (1+m)^2 + \frac{(N+2m)(N-4-2m)}{8} \biggr )
\sum_{i=1}^{k-1} \int_{\Omega} \frac{u^2}{|x|^{2m+4}} X_{1}^{2}
\cdots
X_{i}^{2}\, dx  = \\
= c_N \int_{0}^{1} r^{-1+2\epsilon} X_{1}^{-1+a_1} \cdots
X_{k}^{-1+a_k} \biggl [ \epsilon^2 (2++2m+\epsilon)^2 -
\frac{(N+2m)(N-4-2m)}{2} \epsilon (2+2m+\epsilon) \\
+2 (1+m+\epsilon) \biggl ( - \frac{(N+2m)(N-4-2m)}{4} + \epsilon (
2+2m+ \epsilon) \biggr ) \eta + (1+m+\epsilon)^2 \eta^2 \\
+ 2 \biggl ( - \frac{(N+2m)(N-4-2m)}{4} +\epsilon ( 2+2m+
\epsilon) \biggr ) (\frac{1}{4}\eta^2 +\frac{1}{2}B) \\
- \biggl( (1+m)^2 + \frac{(N+2m)(N-4-2m)}{8} \biggr )
\sum_{i=1}^{k-1} X_{1}^{2} \cdots X_{i}^{2} \biggr ]\, \phi^2\, dr
+O(1),
\end{array}
\]
Using now the identities (\ref{eqb1}), (\ref{eqb2}) and passing to
the limit\ $\epsilon \to 0$, we conclude that
\[
\begin{array}{ccll}
 \int_{\Omega} \frac{|\Delta u|^2}{|x|^{2m}}\, dx - \biggl (
\frac{(N+2m)(N-4-2m)}{4} \biggr )^2 \int_{\Omega}
\frac{u^2}{|x|^{2m+4}} dx \;\;\;\;\;\;\;\;\;\;\;\;\;\;\;
\;\;\;\;\;\;\;\;\;\;\;\;\;\;\;\;\;\;\;\;\;\;\;\;\;\;\;\;\;\;
\;\;\;\;\;\;\;\;\;\;\;\;\;\;\; \\
- \biggl ((1+m)^2 + \frac{(N+2m)(N-4-2m)}{8}\biggr )
\sum_{i=1}^{k-1} \int_{\Omega} \frac{u^2}{|x|^{2m+4}} X_{1}^{2}
\cdots
X_{i}^{2}\, dx  = \\
= c_N ((1+m)^2 + \frac{(N+2m)(N-4-2m)}{8}) \biggl( A_k -
\sum_{i=1}^{k} a_i A_i + \sum_{i=1}^{k-1}\sum_{j=i+1}^{k}
(1-a_j)\Gamma_{ij} \biggr) +O(1).
\end{array}
\]
However, we can pass to the limit\ $a_1 \downarrow 0,...a_{k-1}
\downarrow 0$\ see (\ref{AG}), to conclude that the Rayleigh
quotient now of (\ref{eqbestm1}) is smaller or equal than
\[
\frac{\biggl ( (1+m)^2 + \frac{(N+2m)(N-4-2m)}{8} \biggr ) \biggl(
A_k - a_k A_k \biggr ) + O(1)}{A_k} \to (1+m)^2 +
\frac{(N+2m)(N-4-2m)}{8},
\]
since\ $A_k \to \infty$,\ as\ $a_k \downarrow 0$.\ $\blacksquare$
\vspace{0.2cm} \\
\emph{Proof of Theorem \ref{a.8}}  Is an immediate consequence of
the previous Theorem.
%
%
%
\subsection{The Inequality (\ref{rel2m})}
\setcounter{equation}{0}
In this section we consider Inequality (\ref{rel2m}). For our
approach we consider decomposition into spherical harmonics, see
Section 2. Let\ $u \in C_{0}^{\infty} (\Omega)$.\ Setting\ $u = u
= \sum_{k=0}^{\infty} u_k := \sum_{k=0}^{\infty} f_k (r) \phi_k
(\sigma)$,\ using equalities (\ref{2.3}), (\ref{2.4}) we have that
\begin{eqnarray} \label{9.1}
\int_{\mathbb{R}^N} \frac{|\Delta u_k|^2}{|x|^{2m}}\, dx =
\int_{\mathbb{R}^N} r^{-2m} (f''_{k})^2\, dx + \biggl [
(N-1)(2m+1) +2 c_k \biggr ] \int_{\mathbb{R}^N}
r^{-2-2m} (f'_{k})^2\, dx \nonumber \\
+ c_k \biggl [ c_k + (N-4-2m)(2m+2) \biggr ] \int_{\mathbb{R}^N}
r^{-4-2m} (f_{k})^2\, dx,
\end{eqnarray}
\begin{equation} \label{9.2}
\int_{\mathbb{R}^N} \frac{|\nabla u|^2}{|x|^{2m+2}}\, dx =
\int_{\mathbb{R}^N} r^{-2-2m} (f'_{k})^2\, dx + c_k
\int_{\mathbb{R}^N} r^{-4-2m} (f_{k})^2\, dx.
\;\;\;\;\;\;\;\;\;\;\;\;\;\;\;\;\;\;\;\;\;\;\;\;\;\;\;\;\;\;
\;\;\;\;\;\;\;\;\;
\end{equation}
\setcounter{theorem}{5}
\begin{theorem} \label{p9.1}
Suppose\ $N \geq 5$\ and\ $0 \leq m<\frac{N-4}{2}$.\ Then, for
any\ $u \in C_{0}^{\infty} (\Omega)$,\ the following inequality
holds.
\[
\int_{\Omega} \frac{|\Delta u|^2}{|x|^{2m}}\, dx \geq a_{m,N}\,
\int_{\Omega} \frac{|\nabla u|^2}{|x|^{2m+2}}\, dx,
\]
where\ $a_{m,N}$\ is defined by:
\begin{equation} \label{e9.30}
a_{m,N} := \min_{k = 0,1,2,...} \frac{\biggl (
\frac{(N-4-2m)(N+2m)}{4}+k(N+k-2) \biggr )^2}{\biggl (
\frac{N-4-2m}{2} \biggr )^2 +k(N+k-2)},
\end{equation}
In particular, we have
\[
a_{m,N} = \biggl ( \frac{N+2m}{2} \biggr )^2,
\]
when\ $0 \leq m \leq \frac{-(N+4)+2\sqrt{N^2 -N +1}}{6}$.\
Whereas, we have
\[
a_{m,N} < \biggl ( \frac{N+2m}{2} \biggr )^2,
\]
when\ $\frac{-(N+4)+2\sqrt{N^2 -N +1}}{6} < m < \frac{N-4}{2}$.\
Moreover, the minimum of (\ref{e9.30}) depends only on these\ $k$\
that satisfy
\begin{equation} \label{e9.31}
k \leq \biggl (\frac{\sqrt{3}}{3} - \frac{1}{2} \biggr ) (N-2).
\end{equation}
and let\ $\bar{k}$\ be the largest\ $k$ of (\ref{e9.31}). In
particular, for\ $N \leq 8$,\ and\ $\frac{-(N+4)+2\sqrt{N^2 -N
+1}}{6} < m < \frac{N-4}{2}$\ we have
\[
a_{m,N} = \frac{\biggl ( \frac{(N-4-2m)(N+2m)}{4}+N-1 \biggr
)^2}{\biggl ( \frac{N-4-2m}{2} \biggr )^2 + N-1}
\]
whereas, $8 < N$,\ the interval\ $(\frac{-(N+4)+2\sqrt{N^2 -N
+1}}{6}, \frac{N-4}{2})$\ is been divided in\ $2\bar{k}-1$\
subintervals. For\ $k=1,2,...,\bar{k}$
\begin{eqnarray*}
m^{1}_{k} := \frac{2(N-5)- \sqrt{(N-2)^2-12k(k+N-2)}}{6}, \\
m^{2}_{k} := \frac{2(N-5)- \sqrt{(N-2)^2+12k(k+N-2)}}{6}.
\end{eqnarray*}
When\ $m \in (\frac{-(N+4)+2\sqrt{N^2 -N +1}}{6}, m^{1}_{1}] \cup
[m^{2}_{1}, \frac{N-4}{2})$,\ then
\[
a_{m,N} = \frac{\biggl ( \frac{(N-4-2m)(N+2m)}{4}+N-1 \biggr
)^2}{\biggl ( \frac{N-4-2m}{2} \biggr )^2 + N-1}
\]
For\ $2\leq k \leq \bar{k}-1$\ and\ $m \in (m^{1}_{k},m^{1}_{k+1}]
\cup [m^{2}_{k+1},m^{2}_{k})$,\ then\
\begin{eqnarray*}
a_{m,N} = \min \biggl \{ \frac{\biggl (
\frac{(N-4-2m)(N+2m)}{4}+k(N+k-2) \biggr )^2}{\biggl (
\frac{N-4-2m}{2} \biggr )^2 +k(N+k-2)},\;\;\;\;\;\;\;\;\;\;\;\;\;
\;\;\;\;\;\;\;\;\;\;\;\;\;\;\;\;\;\;\;\;\;\;\;\;\;\;\;\;\;\;\;\;\;
\;\; \\
\;\;\;\;\;\;\;\;\;\;\;\;\frac{\biggl (
\frac{(N-4-2m)(N+2m)}{4}+(k+1)(N+k-1) \biggr )^2}{\biggl (
\frac{N-4-2m}{2} \biggr )^2 +(k+1)(N+k-1)} \biggr \}
\end{eqnarray*}
For\ $m \in (m^{1}_{\bar{k}},m^{2}_{\bar{k}})$,\ then
\begin{eqnarray*}
a_{m,N} = \min \biggl \{ \frac{\biggl (
\frac{(N-4-2m)(N+2m)}{4}+\bar{k}(N+\bar{k}-2) \biggr )^2}{\biggl (
\frac{N-4-2m}{2} \biggr )^2
+\bar{k}(N+\bar{k}-2)},\;\;\;\;\;\;\;\;\;\;\;\;\;
\;\;\;\;\;\;\;\;\;\;\;\;\;\;\;\;\;\;\;\;\;\;\;\;\;\;\;\;\;\;\;\;\;
\;\; \\
\;\;\;\;\;\;\;\;\;\;\;\;\frac{\biggl (
\frac{(N-4-2m)(N+2m)}{4}+(\bar{k}+1)(N+\bar{k}-1) \biggr
)^2}{\biggl ( \frac{N-4-2m}{2} \biggr )^2
+(\bar{k}+1)(N+\bar{k}-1)} \biggr \}
\end{eqnarray*}
Moreover, the constant\ $a_{m,N}$\ in (\ref{rel2m}) is the best.
\end{theorem}
\emph{Proof}\, Decomposing\ $u$\ into spherical harmonics. Using
relations (\ref{9.1}), (\ref{9.2}), and the following Hardy
inequality
\[
\int_{0}^{\infty} r^{N-1-2m} (f''_{k})^2\, dr \geq \biggl (
\frac{N-2-2m}{2} \biggr )^{2} \int_{0}^{\infty} r^{N-3-2m}
(f'_{k})^2\, dr,
\]
we obtain that
\[
a_{m,N} \leq \frac{ C_1 \frac{\int_{0}^{\infty} r^{N-3-2m}
(f'_{k})^2\, dr}{\int_{0}^{\infty} r^{N-5-2m} (f_{k})^2\, dr} +
C_2 }
{\frac{\int_{0}^{\infty} r^{N-3-2m} (f'_{k})^2\,
dr}{\int_{0}^{\infty} r^{N-5-2m} (f_{k})^2\, dr} + c_k}
\]
where\ $C_1 = \biggl [ \biggl ( \frac{N+2m}{2} \biggr )^2 +2 c_k
\biggr ]$,\ $C_2 = c_k \biggl [ c_k - (N-3-2m)(N-4-2m) +
(N-1)(N-4-2m) \biggr ]$.\ However, since\ $C_2 - c_k C_1 \leq 0$,\
the real function
\[
\omega(y) := \frac{C_1 y + C_2}{y + c_k} = \frac{C_1 (y+c_k) + C_2
- c_k C_1 }{y + c_k}
\]
is increasing for positive\ $y$.\ Hence, from the Hardy inequality
\[
\int_{0}^{\infty} r^{N-3-2m} (f'_{k})^2\, dr \geq \biggl (
\frac{N-4-2m}{2} \biggr )^{2} \int_{0}^{\infty} r^{N-5-2m}
(f'_{k})^2\, dr,
\]
we conclude that
\[
a_{m,N} \leq A (k,N,m) := \frac{\biggl (
\frac{(N-4-2m)(N+2m)}{4}+c_k \biggr )^2}{\biggl ( \frac{N-4-2m}{2}
\biggr )^2 +c_k}.
\]
We study the monotonicity of the function
\[
f(x)= \frac{\biggl ( \frac{(N-4-2m)(N+2m)}{4}+ x \biggr
)^2}{\biggl ( \frac{N-4-2m}{2} \biggr )^2 + x},\;\;\; x \geq 0.
\]
It is clear that f admits a (possibly positive) minimum
\[
x_0 := \frac{(N-4-2m)(-N+6m+8)}{4}.
\]

Let\ $N \leq 8$,\ then\ $x_0 \geq 0$.\ In this case holds that\
$c_1 > x_0$,\ for\ $N = 5,6,7,8$\ and\ $0 \leq m \leq
\frac{N-4}{2}$,\ hence\ $a_{m,N} = \min\{ A(0,N,m), A(1,N,m) \}$.\
Comparing\ $A(0,N,m), A(1,N,m)$\ i.e.,\ $A(0,N,m) \leq A(1,N,m)$,\
we obtain that
\begin{equation} \label{m*}
(N+2m)(N-6m-8)+4(N-1) \geq 0.
\end{equation}
By simple calculations we may prove that
\begin{eqnarray*}
&&\mbox{if}\;\; 0 \leq m \leq m^* := \frac{-(N+4)+2\sqrt{N^2 -N
+1}}{6}\;\;\;
\mbox{then}\;\; a_{m,N} = A(0,N,m), \vspace{0.1cm} \\
&&\mbox{while for}\;\; m > m^*,\;\; \mbox{we have that}\;\;
a_{m,N} = A(1,N,m),
\end{eqnarray*}
which in particular holds for every\ $N$.\ In the case where\ $m
\leq \frac{N-8}{6}$,\ (clearly\ $N > 8$),\ we have that\ $x_0 \leq
0$\ and\ $f$\ is increasing for all nonegative\ $x$.\ Hence
\[
a_{m,N} = A(0,N,m),\;\; \mbox{for}\; N >8\;\; \mbox{and}\;\; 0
\leq m \leq \frac{N-8}{6}.
\]
Note that\ $\frac{N-8}{6} < m^*$,\ for every\ $N$.\ In the case
where\ $0< \frac{N-8}{6} \leq m \leq \frac{N-4}{2}$,\ the
situation seems to be more complicated since\ $a_{m,N}$\ also
depends on some\ $k>1$.\ Observe that\ $x_0 > 0$,\ which implies
that\ $f$\ is decreasing for\ $x \in [0, x_0)$\ and increasing
for\ $x_0 < x$.\ In order to estimate the minimum\ $A(k,N,m)$,\ in
terms of\ $k$,\ it suffices to find the relative position of\
$x_0$,\ as\ $c_k$\ varies; Let
\[
\bar{k} := \max\{ k \in \mathbb{N},\;\; \mbox{such that}\; c_k <
x_0 \},
\]
(i.e.\ $c_{\bar{k}} < x_0 < c_{\bar{k}+1})$\ then\ $a_{m,N} = \min
\{ A(\bar{k},N,m), A(\bar{k}+1,N,m) \}$.\ However,\ $c_k < x_0$\
implies that
\[
12m^2 - 8(N-5)m + N^2 -12N +32 +4c_k <0.
\]
Let\ $D := (N-2)^2 -12 c_k$,\ $m^{1}_{k} :=
\frac{2(N-5)-\sqrt{D}}{6}$\ and\ $m^{2}_{k} :=
\frac{2(N-5)+\sqrt{D}}{6}$.\ Then, for every\ $k \in \mathbb{N}$,\
such that\ $D>0$,\ (note that\ $D \ne 0$,\ for any\ $k$,\ $N$)\
there exist a whole interval of\ $m \in (\frac{N-8}{6},
\frac{N-4}{2})$,\ such that\ $a_{m,N} = A(k,N,m)$,\ as follows:
\begin{eqnarray*}
&&\mbox{if}\;\; m \in (m^{1}_{k},m^{1}_{k+1}] \cup
[m^{2}_{k+1},m^{2}_{k}),\;\;\; \mbox{then}\;\; a_{m,N} = \min
\{ A(k,N,m), A(k+1,N,m) \},  \\
\\
&&\mbox{while for}\;\; m \in (m^{1}_{\bar{k}},m^{2}_{\bar{k}}),
\;\;\; \mbox{we have that}\;\; a_{m,N} = \min \{ A(\bar{k},N,m),
A(\bar{k}+1,N,m) \}.
\end{eqnarray*}
Having in mind that\ $m^{1}_{0} = \frac{N-8}{6}$\ and\ $m^{2}_{0}
= \frac{N-4}{2}$,\ we conclude that\ $a_{m,N}$\ behaves in the way
that the theorem states.

Finally, we prove that\ $a_{m,N}$\ is the best constant. To this,
let\ $k$\ be such that\
\[
a_{m,N} =  \frac{\biggl ( \frac{(N-4-2m)(N+2m)}{4}+k(N+k-2) \biggr
)^2}{\biggl ( \frac{N-4-2m}{2} \biggr )^2 +k(N+k-2)}.
\]
We then set
\[
u = |x|^{-\frac{N-4}{2}+m + \epsilon} \phi_k (\sigma) \phi(r),
\]
where\ $\phi(r) \in C^{\infty}_{0} (B_1)$\ is a smooth cutoff
function, such that\ $0 \leq \phi \leq 1$, with\ $\phi \equiv 1$\
in\ $B_{1/2}$\ and\ $\phi_k (\sigma)$\ is an eigenfunction of the
Laplace-Beltrami operator with corresponding eigenvalue\ $c_k =
k(N+k-2)$.\ Then we have that
\[
\frac{1}{c_N} \int_{\Omega} \frac{|\Delta u|^2}{|x|^{2m}}\, dx =
\biggl ( -\frac{(N+2m)(N-4-2m)}{4} - c_k + \epsilon
(2+2m+\epsilon) \biggr )^2 \int_{0}^{1} r^{-1+2\epsilon} \phi^2(r)
dr + O(1),
\]
\[
\frac{1}{c_N} \int_{\Omega} \frac{|\nabla u|^2}{|x|^{2m+2}}\, dx =
\biggl [ \biggl (  -\frac{N-4-2m}{2} + \epsilon \biggr )^2 + c_k
\biggr ] \int_{0}^{1} r^{-1+2\epsilon} \phi^2(r) dr + O(1).
\]
Letting now\ $\epsilon \downarrow 0$\ we obtain the result.\
$\blacksquare$ \vspace{0.2cm}

The requirement (\ref{e9.31}) implies that
\begin{equation} \label{9.3}
k < \frac{2\sqrt{3}-3}{6}\, (N-2),\;\;\; \mbox{or}\;\; k < 0.077\,
(N-2).\
\end{equation}
From (\ref{9.3}) it is clear that for\ $N <15$\ the quantity\
$a_{m,N}$\ depends only on\ $k=0,1$.\ Thus, the case where\
$N=9,...14$,\ is similar to that of\ $N=5,...,8$.\ However, there
is a qualitative difference between these two cases, so we prefer
to state Proposition \ref{p9.1} in this way. Observe that the
above arguments still hold in the case of\ $m=0$,\ see Proposition
\ref{p3.1}.\ As an example assume that\ $N=30$.\ Then, from
(\ref{9.1}) we deduce that\ $\bar{k}=2$.\ We have also that\
$\frac{N-4}{2} = 13$,\ $m* \simeq 4.17$,\ $m^{1}_{1} \simeq
4.85$,\ $m^{1}_{2} = 7$,\ $m^{2}_{2} \simeq 9.66$,\ $m^{2}_{1}
\simeq 11.81$.\ If we take for instance\ $m=8$\ we have that\ $x_0
= 65$,\ when\ $c_2 = 60$\ and\ $c_3 = 93$.\ Actually, in this case
we have that\ $A(k,N,m)=A(0,30,8)= 529$,\ $A(1,30,8)= 384$,\
$A(2,30,8) \simeq 360.29$,\ $A(3,30,8) \simeq 366.64$,\ hence\
$A(30,8)= A(2,30,8)$.\ \vspace{0.3cm} \\
\emph{Proof of Theorem \ref{t9.2}}\, Let\ $V(x) =
\sum_{i=1}^{\infty} X_{1}^{2}(\frac{|x|}{D})
X_{2}^{2}(\frac{|x|}{D}) \ldots X_{i}^{2}(\frac{|x|}{D})$.\ From
relations (\ref{9.1}), (\ref{9.2}), inequality (\ref{9.40}) is
equivalent to
\begin{eqnarray} \label{9.41}
\int_{\Omega} r^{-2m} (f''_{k})^2\, dx - \biggl ( \frac{N-2-2m}{2}
\biggr )^2 \int_{\Omega} r^{-2-2m} (f'_{k})^2\, dx - \frac{1}{4}
\int_{\Omega} r^{-2-2m} V(x) (f'_{k})^2\, dx \nonumber \\
+ 2c_k \int_{\Omega} r^{-2-2m} (f'_{k})^2 dx +c_k \biggl [ c_k +
(N-4-2m)(2m+2) - \biggl ( \frac{N+2m}{2} \biggr )^2 \biggr ]
\int_{\Omega} r^{-4-2m} (f_{k})^2
dx \nonumber \\
- \frac{c_k}{4} \int_{\Omega} r^{-4-2m} V(x) (f_{k})^2\, dx \geq
0.
\end{eqnarray}
However, inequality (\ref{1.4}) implies that
\[
\int_{\Omega} r^{-2m} (f''_{k})^2\, dx - \biggl ( \frac{N-2-2m}{2}
\biggr )^2 \int_{\Omega} r^{-2-2m} (f'_{k})^2\, dx - \frac{1}{4}
\int_{\Omega} r^{-2-2m} V(x) (f'_{k})^2\, dx \geq 0.
\]
Hence, it suffices to prove that
\begin{eqnarray} \label{9.42}
2c_k \int_{\Omega} r^{-2-2m} (f'_{k})^2\, dx +c_k \biggl [ c_k +
(N-4-2m)(2m+2) - \biggl ( \frac{N+2m}{2} \biggr )^2 \biggr ]
\int_{\Omega} r^{-4-2m} (f_{k})^2\,
dx \nonumber \\
- \frac{1}{4} c_k \int_{\Omega} r^{-4-2m} V(x) (f_{k})^2\, dx \geq
0,
\end{eqnarray}
or, since (\ref{9.42}) holds for\ $k=0$,
\begin{eqnarray} \label{9.43}
2 \int_{\Omega} r^{-2-2m} (f'_{k})^2\, dx + \biggl [ c_k +
(N-4-2m)(2m+2) - \biggl ( \frac{N+2m}{2} \biggr )^2 \biggr ]
\int_{\Omega} r^{-4-2m} (f_{k})^2\,
dx \nonumber \\
- \frac{1}{4} \int_{\Omega} r^{-4-2m} V(x) (f_{k})^2\, dx \geq 0,
\end{eqnarray}
for any\ $k=1,2,...$.\ Recalling again inequality (\ref{1.4}),
which gives
\[
\int_{\Omega} r^{-2-2m} (f'_{k})^2\, dx \geq \biggl (
\frac{N-4-2m}{2} \biggr )^2 \int_{\Omega} r^{-4-2m} (f_{k})^2\, dx
+ \frac{1}{4} \int_{\Omega} r^{-4-2m} V(x) (f_{k})^2\, dx,
\]
we obtain that (\ref{9.43}) holds if
\[
2 \biggl ( \frac{N-4-2m}{2} \biggr )^2 + c_k + (N-4-2m)(2m+2) -
\biggl ( \frac{N+2m}{2} \biggr )^2 \geq 0,
\]
for any\ $k=1,2,...$.\ However, this last inequality for\ $k=1$\
is equivalent to (\ref{m*}), which holds for\ $0\leq m \leq
\frac{-(N+4)+2\sqrt{N^2 -N +1}}{6}$.\

Assume now the minimizing sequences
\[
u(x) := w(x)\, \phi(|x|),\;\;\;
w(x) := |x|^{-\frac{N-4}{2}+\epsilon} X_{1}^{\frac{-1+a_1}{2}}
X_{2}^{\frac{-1+a_2}{2}} \cdots X_{k}^{\frac{-1+a_k}{2}},
\]
introduced in Section 5 and using the same notation we have that
\[
\begin{array}{ccll}
\int_{\Omega} \frac{|\Delta u|^2}{|x|^{2m}} dx &=&
\int_{\Omega} \frac{w^2}{|x|^{2m+4}} \biggl [ \biggl ( -
\frac{(N+2m)(N-4-2m)}{4} + \epsilon (2+2m+ \epsilon) \biggr )^2 +
(1+\epsilon+m)^2 \eta^2 \\
&& + 2 (1+m+\epsilon) \biggl ( - \frac{(N+2m)(N-4-2m)}{4}
+\epsilon (2+2m+ \epsilon) \biggr ) \eta \\
&& + 2 \biggl (-\frac{(N+2m)(N-4-2m)}{4}+\epsilon (2+2m+ \epsilon)
\biggr ) (\frac{1}{4} \eta^2 + \frac{1}{2}B) \biggr ] \cdot \cdot
\phi^2\, dx +O(1),
\end{array}
\]
and
\[
\begin{array}{ccll}
\int_{\Omega} \frac{|\nabla u|^2}{|x|^{2m+2}} X_{1}^{2} \cdots
X_{i}^{2}\, dx = \int_{\Omega} \frac{w^2}{|x|^{2m+4}} \biggl [
\biggl ( -\frac{N-4-2m}{2} + \epsilon \biggr )^2 + \biggl (
-\frac{N-4-2m}{2} + \epsilon \biggr ) \eta + \frac{1}{4} \eta^2
\biggr ] \cdot \\
\;\;\;\;\;\;\;\;\;\;\;\;\;\;\;\;\;\;\;\;\;\;\;\;\;\;\;\;\;
\;\;\;\;\;\;\;\;\;\;\;\;\;\;\;\;\;\;\;\;\;\;\;\;\;\;\;\;\; \cdot
X_{1}^{2} \cdots X_{i}^{2}\, \phi^2\, dx +O(1),
\end{array}
\]
We now use identities (\ref{eqb1}), (\ref{eqb2}) and passing to
the limit\ $\epsilon \to 0$, to conclude that
\[
\begin{array}{ccll}
\int_{\Omega} \frac{|\Delta u|^2}{|x|^{2m}} dx - \biggl (
\frac{N+2m}{2} \biggr )^2\, \int_{\Omega} \frac{|\nabla
u|^2}{|x|^{2m+2}} dx - \frac{1}{4} \sum_{i=1}^{k-1} \int_{\Omega}
\frac{|\nabla u|^2}{|x|^{2m+2}} X_{1}^{2} \cdots X_{i}^{2}\, dx
=\;\;\;\;\;\;\;\;\;\;\;\;\;\;\;\;\;\;\;\;\;\;\;\;\;\;
\\
= - \frac{1}{4} \biggl ( \frac{N-4-2m}{2} \biggr )^2 c_N
\int_{0}^{1} r^{-1} X_{1}^{-1+a_1} \cdots X_{i}^{-1+a_k} \biggl [
B + \sum_{i=1}^{k-1} X_{1}^{2} \cdots X_{i}^{2} \biggr ]\,
\phi^2\, dr
+O(1), \\
= \frac{1}{4} \biggl ( \frac{N-4-2m}{2} \biggr )^2 c_N A_k -
\frac{1}{4} \biggl ( \frac{N-4-2m}{2} \biggr )^2 c_N \biggl(
\sum_{i=1}^{k} a_i A_i - \sum_{i=1}^{k-1}\sum_{j=i+1}^{k}
(1-a_j)\Gamma_{ij} \biggr) +O(1).
\end{array}
\]
However, we can pass to the limit\ $a_1 \downarrow 0,...a_{k-1}
\downarrow 0$\ see (\ref{AG}), to conclude that the Rayleigh
quotient now of (\ref{e7.7}) is smaller or equal than
\[
\frac{ \frac{1}{4} \biggl ( \frac{N-4-2m}{2} \biggr )^2 A_k -
\frac{1}{4} \biggl ( \frac{N-4-2m}{2} \biggr )^2 a_k A_k +O(1)}{
\biggl ( \frac{N-4-2m}{2} \biggr )^2 A_k +O(1)} \to \frac{1}{4},
\]
since\ $A_k \to \infty$,\ as\ $a_k \downarrow 0$.\ $\blacksquare$
\vspace{0.2cm} \\
\emph{Proof of Theorem \ref{tp}} Is an immediate consequence of
the previous Theorem.
%
%
%

%

\begin{thebibliography}{MMM}
%
\bibitem[A]{a} Adimurthi, Hardy-Sobolev inequality in $H^1(\Omega)$ and its
applications, \emph{Commun. Contemp. Math.} 4 (2002), no. 3,
409–434.
%
\bibitem[ACR]{acr} Adimurthi, Chaudhuri, Nirmalendu and
Ramaswamy,Mythily, An improved Hardy-Sobolev inequality and its
application, \emph{Proc. Amer. Math. Soc.} 130 (2002), no. 2,
489–505.
%
\bibitem[BT]{bt} Marino Badiale and Gabriella Tarantello, A
Sobolev-Hardy inequality with applications to a nonlinear elliptic
equation arising in astrophysics \emph{Arch. Ration. Mech. Anal.}
163 (2002), no. 4, 259–293.
%
\bibitem[BV]{bv} H. Brezis and J. L. V\'{a}zquez, Blowup solutions
of some nonlinear elliptic problems, \emph{Revista Mat. Univ.
Complutense Madrid} {\bf 10} (1997), 443-469.
%
\bibitem[BFT]{bft03} G. Barbatis, S. Filippas and A. Tertikas,
Series Expansion for\ $L^p$\ Hardy Inequalities,\ \emph{Indiana
Univ. Math. J.} 52 (2003), no. 1, 171--190
%
\bibitem[DH]{dav98} E. B. Davies and A. M. Hinz, Explicit constants for
Rellich Inequalities in\ $L_{p}(\Omega)$,\ \emph{Math. Z.}, 227,
(1998), 511-523.
%
\bibitem[E]{eil01} S. Eilertsen, On weighted franctional integral inequalities,\
\emph{J. Funct. Anal.}, 185, (2001), 342-366.
%
\bibitem[FT]{ft01} S. Filippas and A. Tertikas, Optimizing
Improved Hardy Inequalities, \emph{J. Funct. Anal.}, 192 (2002),
186-233.
%
\bibitem[GGM]{ggm03} F. Gazzola, H. C. Grunau and E. Mitidieri, Hardy inequalities
with optimal constants and remainder terms,\ \emph{Trans. Amer.
Math. Soc.} 356 (2004), no. 6, 2149--2168.
%
\bibitem[GG]{gg} Gabriele Grillo, Hardy and Rellich-Type Inequalities for metrics
Defined by Vector Fields,\ \emph{Potential Analysis} 18 (2003),
187-217.
%
\bibitem[HN]{hn} Y. Han and P. Niu, Hardy-Sobolev type inequalities
on the H-type group, \emph{Manuscripta Math.} 118 (2005), 235-252.
%
\bibitem[MS]{ms} G. Mancini and K. Sandeep, Cylindrical symmetry
of extremals of a Hardy-Sobolev inequality. \emph{Ann. Mat. Pura
Appl.} (4) 183 (2004), no. 2, 165–172.
%
\bibitem[M1]{maz85} V. G. Maz'ja, Sobolev Spaces, \emph{Springer Verlag}, 1985.
%
\bibitem[M2]{maz99} V. G. Maz'ja, The Wiener test for higher order elliptic equations,
\emph{Duke Math J.}, 115, (3), (2002), 479-512.
%
%
\bibitem[TZ]{tz03} A. Tertikas and N. B. Zographopoulos, Optimizing Improved Hardy
Inequalities for the Biharmonic Operator,\ International
Conference on Differential Equations (Hasselt 2003), 1137-1139,
\emph{World Sci. Publishing, River Edge, NJ}, 2005.
%
\bibitem[V]{v} Nicola Visciglia, A note about the generalized
Hardy-Sobolev inequality with potential in
$L^{p,d}(\mathbb{R}^n)$, \emph{Calc. Var. Partial Differential
Equations} 24 (2005), no. 2, 167–184.
%
\bibitem[Y]{yaf99} D. Yafaev, Sharp constants in the Hardy-Rellich inequalities,
\emph{J. Funct. Anal.}, 168, (1999), 121-144.
%
\end{thebibliography}
\end{document}